\documentclass[11pt]{article}
\usepackage{amsmath, amsthm, amsfonts, amssymb, color}
\usepackage{mathrsfs, bbm}

\newtheorem{thm}{Theorem}[section]
\newtheorem{cor}[thm]{Corollary}
\newtheorem{lem}[thm]{Lemma}
\newtheorem{remark}[thm]{Remark}
\newtheorem{prp}[thm]{Proposition}

\theoremstyle{definition}

\definecolor{wco}{rgb}{0.5,0.2,0.3}

\numberwithin{equation}{section} \theoremstyle{remark}

\def\1{{\mathbbm 1}}
\def\non{\nonumber \\}
\def\wt{\widetilde}

\def\R{\mathbb R}

\def\E{\mathbb E}

\def\L{\mathcal L}
\def\S{\mathcal S}

\def\wh{\widehat}

\def\<{\langle}
\def\>{\rangle}

\def\pf{\noindent{\bf Proof.} }

 \def\beq{\begin{equation}}

 \def\P{\mathbb P} 
  \def\ee{\varepsilon}

\begin{document}
\bibliographystyle{plain}

\title{\Large \bf Laplacian Perturbed  by Non-Local Operators}

\author{
{\bf Jie-Ming Wang}\\
\footnotesize{Department of Mathematics, Beijing Institute of Technology, Beijing 100081, China}\\
\footnotesize{E-mail: wangjm@bit.edu.cn}
}

\date{}

\maketitle

\begin{abstract}
Suppose that $d\ge 1$ and $0<\beta<2$. We establish
the existence and uniqueness of the fundamental solution
$q^b(t, x, y)$ to
the operator $\L^b=\Delta+\S^b$,
where
 $$
 \S^bf(x):=
\int_{\R^d} \left( f(x+z)-f(x)- \nabla f(x) \cdot
z\1_{\{|z|\leq 1\}} \right) \frac{b(x, z)}{|z|^{d+\beta}}dz
$$
and $b(x, z)$ is a bounded measurable function on $\R^d\times \R^d$
with $b(x, z)=b(x, -z)$ for $x, z\in \R^d$. We show that if for each $x\in\R^d,$
$b(x, z) \geq 0$ for a.e. $z\in\R^d$, then
$q^b(t, x, y)$ is a strictly positive continuous function
and it uniquely determines a conservative Feller process $X^b$,
which has strong Feller property.
Furthermore, sharp two-sided estimates on $q^b(t, x, y)$ are derived.
\end{abstract}



\section{Introduction}

It is well known that the fundamental solution of the heat equation for Laplacian operator $\Delta$ is Gaussian kernel.
The study of heat kernel estimates for perturbation of Laplace operator $\Delta$ by gradient operator has a long history
and this subject has been studied in many literatures.
In recent years, the study  for nonlocal operator and the associated discontinuous Markov process
 has attracted a lot of interests and
much progress has been made in this field.
In particular,  for  the operator $\Delta$ with a pure nonlocal part in $\R^d$,
Song and Vondracek \cite{SV} obtained the two sided estimates of transition density of the independent sum
of Brownian motion and symmetric stable process, Chen and Kumagai \cite{CK2} generalized the result of \cite{SV}
and established the two sided heat kernel estimates for symmetric diffusion with jumps in a general  setting.
However, it seems that there has been limited literature on the heat kernel for Laplace operator plus a non-symmetric and nonlocal operator until now.
In this paper, our goal is to consider the operator $\Delta$ under a class of nonlocal perturbations in a non-symmetric setting.

Throughout this paper, let  $d\geq 1$ be an integer and $0< \beta <2$.
Recall that a stochastic process $Z=(Z_t, \P_x, x\in \R^d)$
 is called a (rotationally)
symmetric $\beta$-stable process on $\R^d$ if it is
a L\'evy process having
$$ \E_x \left[ e^{i\xi \cdot (Z_t-Z_0)}\right] =e^{-t |\xi|^\beta}
\qquad \hbox{for every } x , \xi \in \R^d.
$$
Let $\wh f(\xi):=\int_{\R^d} e^{i\xi \cdot x } f(x) dx$
 denote the Fourier transform of a function $f$ on $\R^d$.
The fractional Laplacian $\Delta^{\beta/2}$ on $\R^d$ is defined as
\begin{equation}\label{e:1.2n}
 \Delta^{\beta/2} f(x)= \int_{\R^d} \left( f(x+z)-f(x)- \nabla f(x)
 \cdot  z\1_{\{|z|\leq 1\}}  \right) \frac{{\cal A}(d, -\beta)}{|z|^{d+\beta}}dz
\end{equation}
for   $f\in C_b^2(\R^d)$. Here  ${\cal A}(d, -\beta)$
is the normalizing constant so that
$\widehat {\Delta^{\beta/2} f} (\xi)
= -|\xi|^\beta \widehat f (\xi) $.
Hence $\Delta^{\beta/2}$ is the infinitesimal generator
for the symmetric $\beta$-stable process on $\R^d$.

Let $\overline Z_t$ be a finite range symmetric $\beta$-stable process in $\R^d$ with jumps
of size larger than $1$ removed. It is known that the infinitesimal generator of the process $\overline Z_t$ is
the truncated  operator
\begin{equation}\label{e:trun0}
\overline{\Delta}^{\beta/2}f(x):=\int_{|z|\leq 1} \left( f(x+z)-f(x)- \nabla f(x)
 \cdot  z  \right) \frac{{\cal A}(d, -\beta)}{|z|^{d+\beta}}dz.
\end{equation}

Let $b(x, z)$ be a  real-valued bounded function
on $\R^d\times \R^d$ satisfying
\begin{equation}\label{e:b}
b(x, z)=b(x, -z) \qquad \hbox{for every } x, z\in \R^d.
\end{equation}
This paper is concerned with the existence, uniqueness
and sharp estimates on the ``fundamental solution" of
the following  operator  on $\R^d$,
$$
\L^b f(x)=\Delta f(x) + \S^b f(x),
$$
where
 \begin{equation}\label{e:1.1}
\S^bf(x):=
\int_{\R^d} \left( f(x+z)-f(x)- \nabla f(x) \cdot
z \1_{\{|z|\leq 1\}}  \right) \frac{b(x, z)}{|z|^{d+\beta}}dz, \quad f\in C_b^2(\R^d).
\end{equation}
We point out that if  $b(x, z)$ satisfies condition \eqref{e:b},
the truncation $|z|\leq 1$ in \eqref{e:1.1} can be replaced
by $|z|\leq \lambda$
for any $\lambda >0$; that is, for every $\lambda>0$,
\begin{equation}\label{e:1.4a}
\S^bf(x) =
\int_{\R^d} \left( f(x+z)-f(x)-\langle\nabla f(x),
z \> \1_{\{|z|\leq \lambda\}} \right) \frac{b(x, z)}{|z|^{d+\beta}}dz.
\end{equation}
In fact, under condition \eqref{e:b},
\begin{eqnarray}\label{e:1.5a}
\S^bf(x)&=&  \ {\rm p.v.}
\int_{\R^d} \left( f(x+z)-f(x)\right)
\frac{b(x, z)}{|z|^{d+\beta}}dz \nonumber \\
&:=&
\lim_{\varepsilon\to 0}\int_{\{z\in \R^d:
|z|>\varepsilon\}}\left(f(x+z)-f(x)\right)\frac{b(x,z)}{|z|^{d+\beta}}\,dz.
\end{eqnarray}
The operator $\L^b$ is in general non-symmetric.
Clearly,  $\L^b=\Delta$ when $b \equiv 0.$
 $\L^b=\Delta+ \Delta^{\beta /2}$ when $b\equiv  {\cal A}(d, -\beta)$
and $ \L^b=\Delta+\overline \Delta^{\beta /2}$ when $b(x,z)= {\cal A}(d, -\beta)1_{\{|z|\leq 1\}}(z)$.
It is known that the above two symmetric operators are the infinitesimal generators of the independent sum of the
Brownian process   and the symmetric
$\beta$-stable process (rep. symmetric finite range $\beta$-stable process).
The L\'evy measures of $\Delta^{\beta/2}$ and $\overline{\Delta}^{\beta/2}$ are  symmetric in the variable $z$ and do not depend on $x$,
the perturbed operator $\S^b$ under the condition \eqref{e:b} can be viewed as a nonsymmetric extension of
 $\Delta^{\beta/2}$ and $\overline{\Delta}^{\beta/2}.$

Our motivation for  the operator $\L^b$
comes from a very recent work \cite{CW}. The authors  consider
the existence and uniqueness of  the fundamental solution to the fractional Laplacian operator $\Delta^{\alpha/2}$ perturbed
by  lower order nonlocal operator $\S^b$ defined in \eqref{e:1.4a} (i.e. $\Delta^{\alpha/2}+\S^b$ and $0<\beta<\alpha<2$)
and further derive the two sided heat kernel estimates.
The main method to get the upper bound estimate in \cite{CW} is the iterative  Duhamel's formula.
However, this method can't work as well for the operator $\L^b$ as in \cite{CW}.
The main reason is that the Gaussian kernel $p_0(t,x,y)$ is an exponential function of $|x-y|^2,$
after finite number of iterations of recursive Duhamel's formula, the item $\exp(-c^n|x-y|^2/t)$,
where $0<c<1$ is a constant and $n$ is the number of iterations, will appear in the upper bound estimate
and thus one will lose the exponential term in the last.
So  to derive a sharp upper bound estimate, we need a more delicate estimate.

On the other hand, one of main goal of this paper is to study sharp two sided heat kernel estimates for $\Delta$
under finite range non-local perturbation. To the author's knowledge, even for the symmetric case $\Delta+a \overline \Delta^{\beta/2}, a>0,$
the relevant result is new until now. As we shall see in Theorem \ref{T:1.4} and Corollary \ref{C:1.4} below,
 the two sided heat kernel estimates for such operator under some positivity condition will depend
heavily  on the transition density function  for the truncated $\beta/2$-fractional operator $\overline\Delta^{\beta/2}$,
in addition to the Gaussian kernel.
This is different from  the fractional Laplacian operator $\Delta^{\alpha/2}$ under finite range nonlocal  perturbation,
because Theorem 1.4 of \cite{CW} shows that the heat kernel for $\Delta^{\alpha/2}+\S^b$ in this case is comparable with
the heat kernel for fractional Laplacian operator $\Delta^{\alpha/2}$.
Furthermore, the upper bound estimate for $\Delta$
under finite range non-local perturbation in Theorem \ref{T:1.4}  in this paper can't be obtained by the  methods in \cite{CW} and  related literatures,
we use a new probability argument to get it.

For $a\geq 0$, denote by   ${p}_a(t,x, y)$ the
fundamental function   of   $\Delta+a\Delta^{\beta/2}$
.  Clearly,  $  p_a(t, x, y)$
is a function of $t$ and $x-y$, so sometimes we
also  write it as $p_a (t, x-y)$.
It is known (see \eqref{e:h} of  Section 2 for details)
 that on $(0, \infty)\times \R^d \times \R^d$,
\begin{equation}\label{e:1.5}
p_0 (t, x, y)  = (4\pi t)^{-d/2}e^{-|x-y|^2/4t} ,
\end{equation}
\begin{equation}\label{e:1.6t}
\begin{aligned}
&c_1 \left(t^{-d/2}\wedge (at)^{-d/\beta}\right) \wedge \left( p_0(t,c_2x, c_2y)+ \frac{at}{|x-y|^{d+\beta}}\right)\\
\leq & p_a(t, x, y) \leq c_3\left(t^{-d/2}\wedge (at)^{-d/\beta}\right)
 \wedge \left( p_0(t,c_4x, c_4y)+ \frac{at}{|x-y|^{d+\beta}}\right).
 \end{aligned}
\end{equation}
Here we use $a\vee c$ and $a\wedge c$ to denote $\max\{a, c\}$ and $\min \{a, c\}$, respectively.
Note that $(at)^{-d/\beta} \geq t^{-d/2}$
whenever $0<t\leq a^{-2/(2-\beta)}$.
\medskip

To establish the fundamental solution of $\L^b,$ we use the method of Duhamel's formula.
Since $\L^b=\Delta+\S^b$ is a lower order perturbation of $\Delta$ by $\S^b$, the fundamental solution
(or kernel) $q^b(t, x, y)$ of $\L^b$ should satisfy the following
Duhamel's formula:
\begin{equation}\label{e:1.7n}
q^b(t, x, y)=p_0(t, x, y)+\int_0^t \int_{\R^d} q^b (t-s, x, z)
\S^b_z p_0 (s, z, y)dz ds
\end{equation}
for $t>0$ and $x, y\in \R^d$.
Here the notation $S^b_z p_0 (s, z, y)$ means the non-local operator $\S^b$
is applied to the function $z\mapsto p_0(s, z, y)$. Similar notation will also be used for other operators, for example, $\Delta_z$.
Applying \eqref{e:1.7n} recursively, it is reasonable to
conjecture that $\sum_{n=0}^\infty q^b_n(t, x, y)$, if convergent,
 is a solution to \eqref{e:1.7n}, where $q^b_0(t, x, y):=p_0(t, x, y)$
 and
\begin{equation}\label{e:qn}
q^b_n (t, x, y):=\int_0^t \int_{\R^d} q^b_{n-1} (t-s, x, z)
\S^b_z p_0(s, z, y) dz ds \quad \hbox{for } n\geq 1.
\end{equation}

The followings are the main results of this paper.
We use the notation $\|b\|_\infty:= \| b\|_{L^\infty}$.

\begin{thm}\label{T:1.1}
There is a constant $A_0=A_0(d, \beta)>0$
so that for every bounded   function $b$ on
$\R^d\times \R^d$ satisfying condition \eqref{e:b}
with $\| b\|_\infty \leq A_0$,
there is a unique  continuous function $q^b(t,x,y)$
 on $(0, \infty)\times \R^d\times \R^d$
  that   satisfies \eqref{e:1.7n}  on
 $(0, \ee]\times \R^d \times \R^d$
 with $|q^b(t, x, y)| \leq c_1   p_1(t, c_2x, c_2y)$ on
 $(0, \ee]\times \R^d \times \R^d$ for some $\ee, c_1, c_2 >0$, and that
\begin{equation}\label{e:1.4}
\int_{\R^d} q^b(t, x, y) q^b(s, y, z) dy =q^b(t+s, x, z)
\quad \hbox{for every }  t, s >0 \hbox{ and } x, z\in \R^d.
\end{equation}
 Moreover, the following holds.
\begin{description}
\item{\rm (i)}
$q^b(t, x, y)=\sum_{n=0}^\infty q^b_n(t, x, y)$
on $(0,  (1\wedge A_0/\|b\|_\infty)^{2/(2 -\beta)}]\times \R^d \times \R^d$, where
$q^b_n(t, x, y)$ is defined by \eqref{e:qn}.

\item{\rm (ii)}
$q^b(t,x,y)$ satisfies the Duhamel's formula \eqref{e:1.7n} for all $t>0$ and $x,y\in\R^d.$
Moreover, $\S^b_x q^b(t, x, y)$ exists pointwise in the sense of \eqref{e:1.5a} and
\begin{equation}\label{e:Du}
q^b(t, x, y)=p_0(t, x, y)+\int_0^t \int_{\R^d} p_0 (t-s, x, z)
\S^b_z q^b (s, z, y)dz ds
\end{equation}
for $t>0$ and $x, y\in \R^d$.

\item{\rm (iii)}
 For each $t>0$ and $x\in \R^d$,
$\int_{\R^d} q^b(t, x, y) dy =1$.

\item{\rm (iv)} For every $f\in C^2_b (\R^d)$,
$$
T^b_t f(x)-f(x)=\int_0^t T^b_s \L^b f(x) ds,
$$
where $T^b_t f(x)=\int_{\R^d} q^b(t, x, y) f(y) dy$.

\item{\rm (v)} Let   $A>0$.
There are  positive constants $C =C (d, \beta,  A)\geq 1$ and $0<\overline C =\overline C (d, \beta,  A)\leq 1$
so that
for any $b$ with $\| b \|_\infty \leq   A$,
\begin{equation}\label{e:1.7}
|q^b(t,x,y)| \leq C e^{Ct} {p}_{\| b\|_\infty}
(t,\overline Cx,\overline Cy)  \quad \hbox{on }
(0, \infty) \times \R^d\times \R^d.
\end{equation}

\end{description}
\end{thm}

\medskip

In general
the kernel $q^b(t, x, y)$ in Theorem \ref{T:1.1} can
be negative. The next theorem gives the necessary and sufficient condition for $q^b(t,x,y)\geq 0.$
\medskip

\begin{thm}\label{T:1.2}  Let $b$ be a bounded   function on
$\R^d\times \R^d$  that satisfies \eqref{e:b} and that
\begin{equation}\label{e:5.1}
x\mapsto b(x, z) \hbox{ is continuous for a.e. } z\in \R^d.
\end{equation}
Then $q^b(t, x, y)\geq 0$ on $(0, \infty)\times \R^d \times \R^d$
if and only if for each $x\in \R^d$,
\begin{equation}\label{e:1.15}
b(x, z)\geq 0 \quad
\hbox{for  a.e. } z\in \R^d.
\end{equation}
\end{thm}

\bigskip

 Let $\overline{p}_\beta(t,x,y)$ be the fundamental solution of the truncated  operator $\overline{\Delta}^{\beta/2}$
 defined in \eqref{e:trun0}.
It is established in \cite{CKK} that
$\overline p_\beta(t, x, y)$ is jointly continuous and
enjoys the following two sided estimates:
\begin{equation}\label{e:trun1}
 \overline p_\beta(t, x, y) \asymp t^{-d/\beta} \wedge
\frac{t}{|x-y|^{d+\beta}}
\end{equation}
for $ t\in (0, 1]$ and  $|x-y|\leq 1$,
  and there are constants $c_k>0$, $k=1, 2, 3, 4$
  so that
\begin{equation}\label{e:trun2}
c_1 \left( \frac{t}{|x-y|}\right)^{c_2 |x-y|}
\leq \overline p_\beta(t, x, y) \leq
c_3 \left( \frac{t}{|x-y|}\right)^{c_4 |x-y|}
\end{equation}
for $ t\in (0, 1]$ and  $|x-y|> 1$.

\medskip

Next theorem drops the assumption \eqref{e:5.1},
gives lower bound estimate on $q^b (t, x, y)$
for $b(x, z)$ satisfying condition \eqref{e:1.15}
and establish the Feller process
for $\L^b$.

\medskip

Define $m_b:=\inf_x  {\rm essinf}_{z}b(x,z), \, M_b:={\rm esssup}_{x,z}b(x,z).$

\begin{thm}\label{T:1.3}
For every $A>0$, there are positive constants
$C_k=C_k(d,\beta, A), k= 1,\cdots, 4$
such that
for any bounded $b$
with \eqref{e:b}, \eqref{e:1.15} and $\| b\|_\infty\leq A$,
\begin{equation}\label{e:1.14'}
C_1p_{m_b}(t,C_2 x, C_2 y)\leq q^b(t,x,y)
\leq C_3 p_{M_b} (t, C_4x, C_4y)  \quad \hbox{for }
t\in (0,1] \hbox{ and }  x,y\in\R^d.
\end{equation}
Furthermore, for each $\lambda>0$ and $\ee>0,$ there are positive constants
$C_k=C_k(d,\beta, \ee, \lambda, A), k= 5, 6$
such that for any bounded  $b$ on $\R^d\times \R^d$
with $\| b\|_\infty\leq A$ satisfying \eqref{e:b}, \eqref{e:1.15} and
\begin{equation}\label{e:1.8}
\inf_{x\in\R^d, |z|\leq \lambda} b(x,z)>\ee,
\end{equation}
we have
\begin{equation}\label{e:1.14}
C_5  \left(t^{-d/2}\wedge (p_{m_b} (t,C_6x,C_6y)+\overline p_\beta (t,C_6x,C_6y))\right)\leq q^b(t,x,y)
\leq C_3 p_{M_b} (t, C_4x, C_4y)
\end{equation}
for
$t\in (0,1]$ and $x,y\in\R^d$.
For each bounded function $b$ with \eqref{e:b} and \eqref{e:1.15},
the kernel $q^b(t, x, y)$ uniquely determines a Feller process $X^b
=(X^b_t, t\geq 0, \P_x, x\in \R^d)$
on the canonical Skorokhod space
${\mathbb D}([0, \infty), \R^d)$
such that
$$ \E_x \left[ f(X^b_t)\right] =\int_{\R^d} q^b (t, x, y) f(y) dy
$$
for every bounded continuous function $f$ on $\R^d$.
The Feller process $X^b$ is conservative
 and has a L\'evy system $(J^b(x, y)dy, t)$,
where
\begin{equation}\label{e:Jb}
 J^b(x, y)=\frac {b(x, y-x)}{|x-y|^{d+\beta}}.
\end{equation}
\end{thm}

\begin{remark} \rm
 The estimates in \eqref{e:1.14'}  are
 sharp in the sense that
$q^b(t, x, y)=p_0(t, x, y)$ when $b\equiv 0$,
and $q^b(t, x, y)=p_1(t, x, y)$ when $b\equiv {\cal A}(d, -\beta)$.
In particular, it follows from \eqref{e:1.14'}
that for every $A\geq 1 $, there are positive constants
$\wt C_1 , \wt C_2 \geq 1$
so  that for any  $b$ on $\R^d\times \R^d$
satisfying \eqref{e:b} with $ 1/A \leq b(x, z)\leq A$ a.e.
\begin{equation}\label{e:1.16}
(1/ \wt C_1 ) \,  p_1 (t, \wt C_2 x, \wt C_2y)\leq q^b(t,x,y)
\leq \wt C_1 \, {p}_1 (t,x/\wt C_2, y/\wt C_2)  \quad \hbox{for }
t\in (0,1] \hbox{ and }  x,y\in\R^d.
\end{equation}

\end{remark}

\bigskip

If  $b$ is a bounded function satisfying
\eqref{e:b} and \eqref{e:1.15}
so that $b(x, z)=0$ for every $x\in \R^d$ and $|z|\geq R$
for some $R>0$;
or, equivalently if  $\L^b=\Delta+\S^b$ is a lower order perturbation of $\Delta$
by finite range non-local operator $\S^b$,
then  we have the following refined upper bound.

\begin{thm}\label{T:1.4}
For every $\lambda>0$ and $M \geq 1$,
 there  are positive constants $C_k=C_k(d,\beta, M, \lambda), k=7,8$
such that
for any bounded $b$ with \eqref{e:b}, \eqref{e:1.15} and
\begin{equation}\label{e:1.10}
 \sup_x b(x,z)\leq M1_{|z|\leq \lambda}(z),
\end{equation}
we have
\begin{equation}\label{e:1.11}
 q^b(t,x,y) \leq C_7 \left[t^{-d/2}\wedge \left(p_0(t, C_8x, C_8 y)+\overline p_\beta(t, C_8x, C_8y)\right)\right]
 \quad \mbox{for} \: t\in (0,1], \,x,y\in\R^d.\\
\end{equation}
\end{thm}

The following follows immediately from Theorem \ref{T:1.4} and \eqref{e:1.14}.
\begin{cor}\label{C:1.4}
For every $\lambda>0$ and $M \geq 1$,, there  are positive constants $c_k=c_k(d,\beta, M, \lambda), k=1,\cdots, 4$
such that for any bounded $b$ with \eqref{e:b} and
\begin{equation}\label{e:1.9}
M^{-1}1_{|z|\leq \lambda}(z)\leq \inf_x b(x,z)\leq \sup_x b(x,z)\leq M1_{|z|\leq \lambda}(z),
\end{equation}
we have
\begin{equation}\label{e:1.21}
\begin{aligned}
&c_1 \left[t^{-d/2}\wedge \left(p_0(t, c_2x, c_2y)+\overline p_\beta(t, c_2x, c_2y)\right)\right]\\
\leq & q^b(t,x,y) \leq c_3 \left[t^{-d/2}\wedge \left(p_0(t, c_4x, c_4y)+\overline p_\beta(t, c_4x, c_4y)\right)\right]\\
\end{aligned}
\end{equation}
for
$t\in (0,1]$ and   $x,y\in\R^d$.
\end{cor}

\begin{remark}\label{R:1.4} \rm

(i) Corollary \ref{C:1.4} reveals the two sided  heat kernel estimates for  Laplacian
$\Delta$ under the finite range non-local perturbation with the condition \eqref{e:1.9}.
This result seems new even in the symmetric case $\Delta+a \overline{\Delta}^{\beta/2}, a>0$.

\smallskip

(ii) It looks difficult when we try to use Duhamel's formula to get the refined upper bound \eqref{e:1.11}
 in Theorem \ref{T:1.4}
under the assumption that $b(x,z)$ satisfies \eqref{e:1.10}.
The main reason is that the heat kernel $\overline p_\beta(t,x,y)$ of the truncated
operator $\overline \Delta^{\beta/2}$ exhibits  the Poisson type form when $|x-y|>1$ and $t\in (0,1]$ (see \eqref{e:trun2}),
which makes big trouble in the iteration of the recursive Duhamel's formula \eqref{e:qn}.
To circumvent this obstacle, we adopt a new probability argument to go through it .

\end{remark}
\medskip

The rest of the paper is organized as follows. In Section
\ref{S:2}, we derive some estimates on
 $\Delta^{\beta/2}_x p_0(t, x, y)$
that will be used in later.  The existence and uniqueness of the fundamental
solution $q^b(t, x, y)$ of $\L^b$ is given in Section \ref{S:3}.
This is done through a series of lemmas and theorems,
which provide more detailed information on $q^b(t, x, y)$
and $q^b_n(t, x, y)$.
Theorem \ref{T:1.1} then follows from these results.
We show in Section \ref{S:4}
that $\{T^b_t; t>0\}$ is a strongly continuous
semigroup in $C_\infty (\R^d)$.
We then apply Hille-Yosida-Ray theorem and Courr\'ege's first theorem
to establish Theorem \ref{T:1.2}.
When $b$ satisfies \eqref{e:b}, \eqref{e:5.1} and \eqref{e:1.15}, $q^b(t, x, y)$ determines a conservative Feller process $X^b$.
In Section 5,  we extend the result to general bounded $b$ that satisfies \eqref{e:b} and
\eqref{e:1.15}
by approximating it
by a sequence of $\{b^{(n)}, n\geq 1\}$
that satisfy \eqref{e:b}, \eqref{e:5.1} and \eqref{e:1.15}.
Finally, the lower bound estimate in Theorem \ref{T:1.3}  and Theorem \ref{T:1.4} are
established by the L\'evy system of $X^b$ and some probability arguments.

Throughout this paper,
we use the capital letters $C_1,C_2, \cdots $
to denote constants in the statement of the results, and their
labeling will be fixed. The lowercase constants $c_1, c_2, \cdots$
will denote generic constants used in the proofs, whose exact values
are not important and can change from one appearance to another.
 We will use ``$:=$" to denote a definition. For a differentiable function
 $f$ on $\R^d$, we use  $\partial_i f$ and $\partial^2_{ij}f$ to
 denote the partial derivatives $\frac{\partial f}{\partial x_i}$
 and $\frac{\partial^2 f}{\partial x_i \partial x_j}$.

\section{Preliminaries}\label{S:2}

Suppose that $Y$ is a Brownian motion, and $Z$ is a
symmetric $\beta$-stable process on $\R^d$ that is independent of $Z$.
  For any $a \ge 0$, we define $Y^a$ by $Y_t^a:=Y_t+ a^{1/\beta} Z_t$. We will call the process $Y^a$ the independent sum of the
Brownian process  $Y$ and the symmetric
$\beta$-stable process  $Z$ with weight $a^{1/\beta}$. The infinitesimal
generator of $Y^a$ is $\Delta+a \Delta^{\beta/2}$.
Let $p_a (t, x, y)$ denote the transition density of $Y^a$ (or
equivalently the heat kernel of $\Delta+ a
\Delta^{\beta/2}$) with respect to the Lebesgue measure on $\R^d$.
 Recently it is proven in \cite{CK2} and \cite{SV} that
\begin{equation}\label{e:1.0}\begin{aligned}
&c_1 \left(t^{-d/2}\wedge t^{-d/\beta}\right) \wedge \left( p_0(t,c_2x, c_2y)+ \frac{t}{|x-y|^{d+\beta}}\right)\\
\leq & p_1(t, x, y) \leq c_3\left(t^{-d/2}\wedge t^{-d/\beta}\right)
 \wedge \left( p_0(t,c_4x, c_4y)+ \frac{t}{|x-y|^{d+\beta}}\right)
 \end{aligned}
 \end{equation}
 for $(t,x,y)\in (0,\infty)\times\R^d\times\R^d.$

Unlike the case of  the Brownian motion $Y:=Y^0$,
$Y^a$ does not have the stable  scaling for $a>0$. Instead, the
following approximate scaling property holds : for every
$\lambda>0$, $\{\lambda^{-1} Y^{a}_{\lambda^2 t}, t\geq 0\}$ has the same distribution as $\{Y^{a \lambda^{(2-\beta)} }_t, t\geq 0\}$.
Consequently, for any $\lambda>0$, we have
\begin{equation}\label{e:scaling}
p_{a\lambda^{(2-\beta) }}  ( t,  x, y) =
\lambda^d p_{a}  (\lambda^{2}t, \lambda x, \lambda y) \qquad
\hbox{for } t>0 \hbox{ and } x, y \in \R^d.
\end{equation}
In particular, letting $a=1$, $\lambda= a^{1/(2 -\beta)},$
 we get
$$
p_a(t,x,y)= a^{d/(2-\beta)}p_1(a^{2/
 (2-\beta)} t ,  a^{1/(2-\beta)} x ,
a^{1/(2-\beta)} y) \qquad \hbox{for } t>0
\hbox{ and } x, y\in \R^d.
$$
So we deduce from \eqref{e:1.0} that there exist
constants $c_k, k=1,\cdots, 4$ depending only on
 $d$ and $\beta$ such that for every $a>0$
 and $(t,x,y) \in (0,
\infty)\times \R^d \times \R^d$
\begin{equation}\label{e:h}
\begin{aligned}
&c_1 \left(t^{-d/2}\wedge (at)^{-d/\beta}\right) \wedge \left( p_0(t,c_2x, c_2y)+ \frac{at}{|x-y|^{d+\beta}}\right)\\
\leq & p_a(t, x, y) \leq c_3\left(t^{-d/2}\wedge (at)^{-d/\beta}\right)
 \wedge \left( p_0(t,c_4x, c_4y)+ \frac{at}{|x-y|^{d+\beta}}\right).
 \end{aligned}
\end{equation}
In fact, \eqref{e:h} also holds when $a=0$.

\bigskip
Recall that $p_0(t, x-y)$ is the transition density function of
Brownian motion $Y$.

\bigskip

\begin{lem}\label{0}
There exists a constant $C_{9}=C_{9}(d)>0$
such that for every $t>0$, $x\in \R^d$ and
$i, j=1, \dots, d$,
$$
p_0(t, x) \leq C_{9} t^{-d/2} \left(1\wedge \frac{t^{1/2}}{|x|}\right)^{d+2}, \quad
| \partial^2_{ij}   p_0 (t,x) |\leq
C_{9} t^{-(d+2)/2} \left(1\wedge \frac{t^{1/2}}{|x|}\right)^{d+4}.
$$
\end{lem}

\pf
It is known that
$$p_0(t,x)=(4\pi t)^{-d/2}e^{-|x|^2/4t}.$$
Thus, $p_0(t,x)\leq (4\pi t)^{-d/2}.$
On the other hand, by the proof of Lemma 2.1 in \cite{SV},
\begin{equation}\label{e:p00}
p_0(t,x)\leq c\frac{t}{|x|^{d+2}}, \quad t>0, \, x\in\R^d.
\end{equation}
For the reader's convenience, we spell out the details here.
For each $x\neq 0,$ define $f: (0,\infty)\rightarrow (0,\infty)$ as follows:
$$f(t)=t^{-1-d/2}\exp \left(-\frac{|x|^2}{4t} \right).$$
Then $f(0+)=f(+\infty)=0.$ Further
$$f'(t)=f(t)t^{-2}(-(d/2+1)t+|x|^2/4).$$
This derivative is zero for $t_0=\frac{|x|^2}{4(d/2+1)},$ positive for $t<t_0$ and negative for $t>t_0.$
Thus, $\max f(t)\leq f(t_0)=c |x|^{-(d+2)}$ and \eqref{e:p00} follows from it.
Therefore, there exists $c_1>0$ such that
\begin{equation}\label{e:p0}
p_0(t,x)\leq c_1\left(t^{-d/2}\wedge \frac{t}{|x|^{d+2}}\right).
\end{equation}

Next, we prove the inequality about the second derivatives of $p_0(t,x).$
By simple computation and \eqref{e:p0}, we have
$$\begin{aligned}\left|\partial^2_{ij}p_0(t,x)\right|
&\leq  \left[\frac{|x|^2}{t^2}+\frac{2}{t}\right]p_0(t,x)\\
&=(4\pi)^2|x|^2 p^{(d+4)}_0(t,\tilde{x}_1)+8\pi p^{(d+2)}_0(t,\tilde{x}_2)\\
&\leq c_2|x|^2\left(t^{-(d+4)/2}\wedge \frac{t}{|x|^{d+6}}\right)+c_2\left(t^{-(d+2)/2}\wedge \frac{t}{|x|^{d+4}}\right)\\
&\leq c_3\left(t^{-(d+2)/2}\wedge \frac{t}{|x|^{d+4}}\right)
\end{aligned}$$
where $\tilde{x}_1\in\R^{d+4}$ and $\tilde{x}_2\in\R^{d+2}$ such that $|\tilde{x}_1|=|\tilde{x}_2|=|x|,$
$p^{(d+4)}_0$ and $p^{(d+2)}_0$ are the transition densities of Laplacian
operator in dimension $d+4$ and $d+2$.
\qed

\bigskip

 Define for $t>0$ and $x\in \R^d$,  the function
$$
|\Delta^{\beta/2}_x| p_0(t, x)
 \begin{cases}
=\int_{|z|\leq t^{1/2}}
\big| p_0 (t, x+z)-p_0(t,x)- \frac{\partial}{\partial x}p_0(t, x)
\cdot z  \big| \, \dfrac{1} {|z|^{d+\beta}}dz \\
 \hskip 0.3truein
  +\int_{|z|> t^{1/2}}|p_0(t,x+z)-p_0(t,x)|\dfrac{dz}{|z|^{d+\beta}}
   \hskip 0.8truein \hbox{for } |x|^2\leq t,\\
= \int_{|z|\leq |x|/2}|p_0(t,x+z)-p_0(t,x)- \frac{\partial}{\partial x}p_0(t,x)\cdot z |\dfrac{1}{|z|^{d+\beta}} dz\\
 \hskip 0.3truein  +\int_{|z|>|x|/2}|p_0(t,x+z)-p_0(t,x)|\dfrac{dz}{|z|^{d+\beta}}
 \hskip 0.8truein \hbox{for } |x|^2> t .
\end{cases}
$$
Let
\begin{equation}\label{e:3.9}
f_0(t, x, y):=\left( t^{1/2} \vee |x-y| \right)^{-(d+\beta)}
= t^{-(d+\beta)/2} \left( 1 \wedge \frac{t^{1/2}}{|x-y|} \right)^{d+\beta}.
\end{equation}

\begin{lem}\label{1}
There exists a constant $C_{10}=C_{10}(d, \beta )>0$ such that
\begin{equation}\label{e:3.2}
 |\Delta^{\beta/2}_x| p_0(t,x,y) \leq C_{10} f_0(t, x, y)
\qquad \hbox{on } (0, \infty)\times \R^d \times \R^d.
\end{equation}
\end{lem}

\pf  (i) We first consider the case $|x|^2\leq t.$ In
this case,
$$\begin{aligned}
|\Delta^{\beta/2}_x| p_0(t,x)&= \int_{|z|\leq t^{1/2}}|p_0(t,x+z)-p_0(t,x)- \frac{\partial}{\partial x}p_0(t,x)
\cdot z |\dfrac{dz}{|z|^{d+\beta}}\\
&\qquad +\int_{|z|\geq t^{1/2}}|p_0(t,x+z)-p_0(t,x)|\dfrac{dz}{|z|^{d+\beta}} \\
&=I+II.
\end{aligned}$$

Note that by Lemma \ref{0},
$$\sup_{u\in \R^d} \left|\frac{\partial^2}{\partial u_i\partial u_j}p_0(t, u)\right|\leq C_9t^{-(d+2)/2}, $$
so  by Taylor's formula,
$$
I \leq  \sup_{u\in \R^d} \left|\frac{\partial^2}{\partial u_i\partial u_j}p_0(t, u)\right|
\int_{|z|\leq t^{1/2}} \frac{|z|^2}{|z|^{d+\beta}}\,dz
\leq c_1t^{-(d+2)/2}t^{(2-\beta)/2}\leq
c_1t^{-(d+\beta)/2}.
$$
For the second item $II,$ we have
$$
II\leq    \int_{|z|\geq t^{1/2}}
\left(p_0(t,x+z)+p_0(t,x)\right) \frac{dz}{|z|^{d+\beta}}
 \leq c_2t^{-d/2}\int_{|z|\geq t^{1/2}}
\frac{1}{|z|^{d+\beta}}dz \leq c_3t^{-(d+\beta)/2}.
$$

(ii) Next, we consider the case $|x|^2\geq t.$ In this case,
$$\begin{aligned}
|\Delta^{\beta/2}_x| p_0(t,x)&=  \int_{|z|\leq |x|/2}|p_0(t,x+z)-p_0(t,x)- \frac{\partial}{\partial x}p_0(t,x)\cdot z |\dfrac{dz}{|z|^{d+\beta}}\\
&\qquad + \int_{|z|\geq |x|/2}|p_0(t,x+z)-p_0(t,x)|\dfrac{dz}{|z|^{d+\beta}}\\
&=:I+II.
\end{aligned}$$
Note that $|x+z|\geq |x|/2$ for $|z|\leq |x|/2$. So by Lemma \ref{0},
$$\sup_{|z|\leq |x|/2}\Big|\frac{\partial^2}{\partial x_i\partial x_j}p_0(t,x+z) \Big|
\leq c_3\sup_{|z|\leq |x|/2}t|x+z|^{-(d+4)}\leq 2^{(d+4)}c_3t|x|^{-(d+4)}.$$ Hence, by
Taylor's formula
\begin{equation}\label{e:w1}
\begin{aligned}
I&\leq  \sup_{|z|\leq |x|/2} \Big|
\frac{\partial^2}{\partial x_i\partial x_j}p_0(t, x+z) \Big|
\int_{|z|\leq |x|/2} \frac{|z|^2}{|z|^{d+\beta}}\,dz\\
&\leq c_4t|x|^{-(d+4)}|x|^{2-\beta}=
c_4t|x|^{-(d+2+\beta)}.
\end{aligned}\end{equation}
As $|x|^2\geq t,$ thus $I\leq c_4|x|^{-(d+\beta)}.$
On the other hand, noting that
$p_0(t,y)\leq p_0(t,x)$ if $|y|\geq |x|.$
Hence, by  the condition that
$|x|^2\geq t,$ we obtain
\begin{equation}\label{e:w2}\begin{aligned}
II&\leq  \int_{|z|\geq |x|/2, |x+z|\geq |x|}
2p_0(t,x)\frac{dz}{|z|^{d+\beta}}
+ \int_{|z|\geq |x|/2, |x+z|\leq |x|} 2p_0(t,x+z)\frac{dz}{|z|^{d+\beta}}\\
&\leq 2 p_0(t,x)\int_{|z|\geq |x|/2}\frac{dz}{|z|^{d+\beta}}
+2^{d+1+\beta} |x|^{-(d+\beta)}\int_{z\in\R^d}p_0(t,x+z)\,dz\\
&\leq c_5t|x|^{-(d+2)}|x|^{-\beta}+2^{d+1+\beta} |x|^{-(d+\beta)}\leq
c_6|x|^{-(d+\beta)}.
\end{aligned}\end{equation}
This establishes the lemma. \qed

\bigskip

\begin{lem}\label{n1}
There is a constant $C_{11}=C_{11}(d, \beta)>0$
such that
\begin{equation}\label{L1}
\int_0^t\int_{\R^d}f_0(s, z, y) dz ds  \leq C_{11} \, t^{1-\beta/2}, \qquad  t\in(0, \infty),\, y\in\R^d.
\end{equation}
\end{lem}

\pf  By the definition of $f_0,$
 \begin{eqnarray*}
&&\int_0^t\int_{\R^d}f_0(s, z, y)\,dz\,ds\\
&\leq& \int_0^t\int_{|y-z|\leq s^{1/2}} s^{-(d+\beta)/2}\,dz\,ds
+\int_0^t\int_{|y-z|> s^{1/2}} \frac{1}{|y-z|^{d+\beta}}\,dz\,ds \\
&\leq & c_1 \int_0^t s^{-\beta/2}\,ds \leq c_2t^{1-\beta/2}.
\end{eqnarray*}
\qed

Define
\begin{equation}\label{e:3.3}
h(t,x,y)=t^{-d/2}\wedge \left(p_0(t,x,y)+\frac{t}{|x-y|^{d+\beta}}\right).
\end{equation}
Then $\int_{\R^d} h(t,x,y)\,dy \asymp 1.$
Moreover, for $t\in (0,1]$, $h(t,x,y)\asymp t^{-d/2}$ when $|x-y|\leq t^{1/2}$ and $h(t,x,y)\asymp p_0(t,x,y)+\frac{t}{|x-y|^{d+\beta}}$
when $|x-y|>t^{1/2}.$ Here for two non-negative functions $f$ and $g$,
the notation $f\asymp g$ means that there is a constant $c\geq 1$ so that $c^{-1} f\leq g\leq cf$ on
their common domain of definitions.

\medskip

\begin{lem}\label{2'}
There exist $C_{12}=C_{12}(d,\beta)>1$ and $0<C_{13}=C_{13}(d,\beta)<1$ such that for any $t\in (0,1],$
$$\int_0^t \int_{\R^d}h(t-s,x,z)f_0(s,z,y)\,dz\,ds
\leq \left\{\begin{array}{ll}
C_{12} h(t,x,y), & \, |x-y|\leq t^{1/2}, \, \mbox{or}\,  |x-y|>1,\\
C_{12} h(t, C_{13}x, C_{13}y), & \, t^{1/2}< |x-y|\leq 1
\end{array}\right. .$$
\end{lem}

\pf  Denote by $I=\int_0^t\int_{\R^d} h(t-s,x,z)
f_0(s,z,y)\,dz\,ds.$

(i) Suppose that $|x-y|\leq t^{1/2}$.
We write $I$ as
 \begin{eqnarray*}
 I&=&\int_0^{t/2}\int_{\R^d}
h(t-s,x,z)f_0(s,z,y)\,dz\,ds\\
&&\quad+\int_{t/2}^{t}\int_{\R^d}
h(t-s,x,z) f_0(s,z,y)\,dz\,ds\\
&=&I_{1}+I_{2}.
\end{eqnarray*}

If $s\in (0,t/2),$ then $t-s\in [t/2,t)$.
Thus
$h(t-s,x,z)\leq c_1t^{-d/2}$.
Hence, by Lemma \ref{n1},
$$
I_{1}\leq c_1t^{-d/2}\int_0^t\int_{\R^d}f_0(s,z,y)\,dz\,ds
\leq c_2  \, t^{-d/2}.
$$
When $s\in [t/2,t]$, noting that $f_0(s,z,y)\leq c_3t^{-(d+\beta)/2}$,
hence,
$$
I_{2}\leq c_3t^{-(d+\beta)/2}\int_0^t\int_{\R^d}h(t-s,x,z)\,dz\,ds
\leq c_4 \, t^{-d/2}.
$$
We thus conclude from the above that there is a $c_5>0$ such that $I\leq c_5 \, h(t,x,y)$
for every $t\in (0,1]$ whenever $|x-y|\leq t^{1/2}$.

(ii) Next assume that $|x-y|\geq t^{1/2}$.  Then
$$\begin{aligned}
I &=\int_0^t\int_{|x-z|\leq |x-y|/2} h(t-s,x,z)f_0(s,z,y)\,dz\,ds\\
&\quad+\int_0^t\int_{|x-z|>|x-y|/2}
h(t-s,x,z)f_0(s,z,y)\,dz\,ds\\
&=:I_1+I_2.
\end{aligned}$$

If $|x-z|\leq |x-y|/2,$ then $|y-z|\geq |x-y|/2>t^{1/2}/2$.
Thus  $f_0(s,z,y)\leq
c_8|x-y|^{-(d+\beta)}$ for $s\in (0,t).$
 Therefore,
 \begin{equation}\label{e:fu1}\begin{aligned}
 I_1&= \int_0^t\int_{|x-z|\leq |x-y|/2}
h(t-s,x,z)f_0(s,z,y)\,dz\,ds\\
& \leq c_6|x-y|^{-(d+\beta)}\int_0^t\int_{\R^d}
h(t-s,x,z)\,dz\,ds\\
&\leq c_6\frac{t}{|x-y|^{d+\beta}}\\
&\leq c_6h(t,x,y)\\
\end{aligned}\end{equation}

Now we consider $I_2$ where $|x-z|>|x-y|/2$. If $|x-y|>1$, then   by Lemma \ref{0},
$p_0(t-s,x,z)\leq c_7\frac{t}{|x-z|^{d+2}}\leq c_7 2^{d+2}\frac{t}{|x-y|^{d+2}}\leq c_7 2^{d+2}\frac{t}{|x-y|^{d+\beta}}$.
So by \eqref{e:3.3}, for $|x-z|>|x-y|/2>1/2,$
$$h(t-s,x,z)\leq \left[ p_0(t-s, x, z)+|x-z|^{-(d+\beta)}t\right]
\leq c_{8}t|x-y|^{-(d+\beta)}\leq c_{8}h(t,x,y).$$
Thus,
$$
I_{2}
\leq c_{8}h(t,x,y)\int_0^t\int_{\R^d}f_0(s,z,y)\,dz\,ds\leq c_{9}h(t,x,y), \quad |x-y|>1.
$$
Therefore, combining the above inequality with \eqref{e:fu1}, there exists $c_{10}>0$ so that
\begin{equation}\label{e:2}
I\leq c_{10}h(t,x,y), \quad |x-y|>1.
\end{equation}

On the other hand, if $t^{1/2}< |x-y|\leq 1,$ then
we divide $I_2$ into two parts:
\begin{equation}\label{e:fu4} \begin{aligned}
I_2
&=\int_0^{t/2}\int_{|x-z|>|x-y|/2} h(t-s,x,z) f_0(s,z,y)\,dz\,ds\\
&\qquad+\int_{t/2}^t\int_{|x-z|>|x-y|/2} h(t-s,x,z) f_0(s,z,y)\,dz\,ds\\
 &=I_{21}+I_{22}.
\end{aligned}\end{equation}
We first consider $I_{21}.$  Noting that $h(t-s, x, z)\leq c_{11} h(t, x/2, y/2)$ for $|x-z|>|x-y|/2$ and $s\in (0, t/2],$ we have
\begin{equation}\label{e:1}\begin{aligned}
I_{21}\leq &c_{11}h(t,x/2,y/2)\int_0^{t/2}\int_{|x-z|>|x-y|/2}  f_0(s,z,y)\,dz\,ds\\
\leq & c_{12} h(t,x/2,y/2), \qquad\qquad t^{1/2}< |x-y|\leq 1.
\end{aligned}\end{equation}
Next, note that by \eqref{e:1.0}, there are constants $c_{13}>1$ and $0<c_{14}<1$ so that
$h(s,x,z)\leq c_{13}p_1(s, c_{14}x, c_{14}z)$ for $s\in (0,1]$ and $x,z\in\R^d.$  Moreover,
$$
f_0(s,z,y)\leq \frac{1}{s}\left[s^{-d/2}\wedge \frac{s}{|y-z|^{d+\beta}}\right]
\leq \frac{1}{s}h(s,z,y)\leq \frac{2c_{13}}{t} p_1(s, c_{14}z, c_{14}y), \quad s\in (t/2, t],
$$
then we have,
\begin{equation}\label{e:1'}\begin{aligned}
I_{22}
&\leq 2c_{13}^2\frac{1}{t}\int_{t/2}^t\int_{\R^d} p_1(t-s, c_{14}x, c_{14}z) p_1(s, c_{14}z, c_{14}y)\,dz\,ds\\
&\leq c_{15} p_1(t, c_{14}x, c_{14}y)\\
&\leq c_{16}  h(t, c_{17}x, c_{17}y), \qquad  t^{1/2}< |x-y|\leq 1,
\end{aligned}\end{equation}
where the constant $c_{17}$ is less than $1$ and the last inequality holds due to \eqref{e:1.0}.
By \eqref{e:fu1} and \eqref{e:fu4}-\eqref{e:1'}, there are $c_{18}>1$ and $0<c_{19}<1$ so that
$$I\leq c_{18} h(t, c_{19}x, c_{19}y), \quad t^{1/2}<|x-y|\leq 1.$$
Therefore, the proof is complete.
 \qed

\section{Fundamental solution}\label{S:3}

Throughout the  rest of this paper, $b(x, z)$ is a bounded function on
$\R^d\times \R^d$ satisfying condition \eqref{e:b}.
Recall the definition of the non-local operator $\S^b$ from \eqref{e:1.1}.
Let $|q^b|_0(t,x,y)=p_0(t,x,y),$
and define for each $n\geq 1$,
$$
|q^b|_n(t,x,y)=\int_0^t\int_{\R^d} |q^b|_{n-1}(t-s,x,z)|\S^b_z p_0
(s,z,y)|\,dzds.
$$
Note that by Lemma \ref{1} and \eqref{e:b},
\begin{equation}\label{L}
|\S^b_z p_0(s,z,y)|\leq \|b\|_\infty |\Delta_z^{\beta/2}|p_0(s,z,y)\leq C_{10}\|b\|_\infty f_0(s,z,y), \quad s>0, \,z,y\in\R^d.
\end{equation}
In view of \eqref{e:3.3}, there is a  constant $C_{14}>1$  so that
\begin{equation}\label{L0}
p_0(t,x,y)\leq C_{14}h(t, x, y), \quad t>0,\, x,y\in\R^d.
\end{equation}

\begin{lem}\label{L:3.3'}
For each $n\geq 0$ and every bounded function $b$ on $\R^d\times\R^d$ satisfying condition \eqref{e:b}, there exists a finite constant $C(n)$ depending on $n$ so that
\begin{equation}\label{k'}
|q^b|_n(t,x,y)
\leq C(n)h(t, C_{13}^n x, C_{13}^n y)<\infty, \quad t\in (0,1],\, x,y\in\R^d.
\end{equation}
\end{lem}

\pf We prove this lemma by induction.
 \eqref{k'} clearly holds for $n=0$ by \eqref{L0}.
Suppose that \eqref{k'} holds for $n=j\geq 0$.  By
 Lemma \ref{2'} and the fact that $0<C_{13}<1,$ we have
 $$\int_0^t\int_{\R^d} h(t-s,  x, z)
f_0 (s, z, y)\,dz\,ds\leq C_{12}h(t, C_{13}x, C_{13}y), \quad t\in (0,1],\, x,y\in\R^d.$$
Then by the above inequality  and (\ref{L}), for $t\in (0,1],\, x,y\in\R^d,$
\begin{equation}\begin{aligned}\label{e:qb1}
& |q^b|_{j+1}(t,x,y)\\
\leq& C(j)
 \,
\int_0^t\int_{\R^d} h(t-s, C_{13}^jx, C_{13}^j z)  |\S^{b}_z p_0(s,z,y)|\,dz\,ds\\
\leq& C(j)
\, C_{10}\|b\|_\infty \,
\int_0^t\int_{\R^d} h(t-s,  C_{13}^jx, C_{13}^j z)
f_0 (s, C_{13}^jz, C_{13}^j y)\,dz\,ds\\
\leq& C(j) \|b\|_\infty C_{10}C_{12}C^{-j}_{13}h(t, C_{13}^{j+1}x, C_{13}^{j+1}y),
\end{aligned}\end{equation}
where the second  inequality holds due to \eqref{L} and $ 0< C_{13}<1.$
Let $C(j+1)=C(j) \|b\|_\infty C_{10}C_{12}C^{-j}_{13},$ then the proof is complete.
\qed

\bigskip

Now we define $q_n^b: (0, 1]\times\R^d\times\R^d\rightarrow \R$ as follows.
 For $t>0$ and $x,y\in\R^d,$
let $q^b_0(t,x,y)=p_0(t,x,y),$ and for each $n\geq 1$, define
\begin{equation}\label{e:3.6}
q_n^b(t,x,y)=\int_0^t\int_{\R^d} q^b_{n-1}(t-s,x,z)\S^b_z p_0(s,z,y)\,dz\,ds.
\end{equation}
Clearly by Lemma \ref{L:3.3'},
each $q_n^b(t,x,y)$ is well defined on $(0, 1]\times\R^d\times\R^d$.

\medskip

For $\lambda >0$, define
\begin{equation}\label{e:4.1}
b^{(\lambda)} (x, z)= \lambda^{\beta/2 -1}
b(\lambda^{-1/2} x, \lambda^{-1/2} z).
\end{equation}
For a function $f$ on $\R^d$, set
$$ f^{(\lambda)} (x):= f(\lambda^{-1/2} x).
$$
By a change of variable, one has from  \eqref{e:1.1}
that
$$ \Delta f^{(\lambda)} (x)= \lambda^{-1} (\Delta  f) (\lambda^{-1/2}x)
$$
and
\begin{equation}\label{e:4.2}
 \S^{b^{(\lambda)}} f^{(\lambda)} (x)= \lambda^{-1} ( \S^b f) (\lambda^{-1/2}x).
\end{equation}
Note that the transition density function $p_0(t, x, y)$ of
the Brownian motion has the following scaling property:
\begin{equation}\label{e:4.3}
p_0(t, x, y)= \lambda^{-d/2} p_0(\lambda^{-1}t, \lambda^{-1/2} x, \lambda^{-1/2} y)
\end{equation}
Recall  $q^b_n(t, x, y)$ is the function defined inductively by \eqref{e:3.6} with $q^b_0(t, x, y):=p_0(t, x, y)$.

\begin{lem}\label{L:3.5b}
Suppose that $b$ is a bounded function on $\R^d\times \R^d$ satisfying \eqref{e:b}. For  every integer $n\geq 0$,
\begin{equation}\label{e:4.4}
 q^{b^{(\lambda)}}_n (t, x, y)=
\lambda^{-d/2} \, q^{b}_n  (\lambda^{-1}t, \lambda^{-1/2} x, \lambda^{-1/2} y), \qquad t\leq 1\wedge\lambda,  x, y\in \R^d;
\end{equation}
or, equivalently,
\begin{equation}\label{e:4.5}
q^b_n (t, x, y)=
\lambda^{d/2} \, q^{b^{(\lambda)}}_n  (\lambda t, \lambda^{1/2} x, \lambda^{1/2} y), \qquad t\leq 1\wedge\lambda^{-1}, x, y\in \R^d.
\end{equation}
\end{lem}

\pf We prove it by induction. Clearly in view of \eqref{e:4.3},
\eqref{e:4.4} holds when $n=0$. Suppose that \eqref{e:4.4} holds
for $n=j\geq 0$. Then by the definition \eqref{e:3.6}, \eqref{e:4.2}
and \eqref{e:4.3},
\begin{eqnarray*}
&& q^{b^{(\lambda)}}_{j+1} (t, x, y) = \int_0^t \int_{\R^d}
q^{b^{(\lambda)}}_j (t-s, x, z) \S^{b^{(\lambda)}}_z p_0(s, z , y)\, dz ds \\
&=&  \int_0^t \int_{\R^d} \lambda^{- d/2}
 q^{b}_j  (\lambda^{-1}(t-s), \lambda^{-1/2} x, \lambda^{-1/2} z)
\lambda^{-d/2-1} \left( \S^b_z p_0 (\lambda^{-1}s, \cdot , \lambda^{-1/2}y)\right) (\lambda^{-1/2}z) \, dz ds \\
&=& \lambda^{-d/2} \int_0^{\lambda^{-1}t} \int_{\R^d}
 q^{b}_j  (\lambda^{-1}t-r, \lambda^{-1/2} x, w)
 \left( \S^b_w p_0 (r, \cdot , \lambda^{-1/2} y)\right)(w) \,  dw dr\\
&=& \lambda^{-d/2}  q^{b}_{j+1} (\lambda^{-1}t, \lambda^{-1/2} x, \lambda^{-1/2} y).
\end{eqnarray*}
This proves that \eqref{e:4.4} holds for $n=j+1$ and so, by induction,
it holds for every $n\geq 0$. \qed

In the following, we use Lemma \ref{L:3.5b} together with Lemma \ref{2'} to get the refined upper bound of
$|q^b_n(t,x,y)|.$

\begin{lem}\label{L:3.3}
For each $A>0$ and every bounded function $b$ on $\R^d\times\R^d$ satisfying condition \eqref{e:b}
with $\|b\|_\infty\leq A,$
\begin{equation}\label{k}
|q^b_n(t,x,y)|
\leq C_{14} \left(AC_{10}C_{12}  \right)^{n}h(t, x, y), \quad t\in (0,1],\, x,y\in\R^d, \, n\geq 0.
\end{equation}
\end{lem}

\pf We prove this lemma by induction.
By \eqref{L0} ,   \eqref{k} clearly holds for $n=0$.
Suppose that \eqref{k} holds for $n=j\geq 0$.  Then by
 Lemma \ref{2'} and (\ref{L}), for $t\in (0,1],\, |x-y|\leq t^{1/2}\,  \mbox{or} \, t\in (0,1], |x-y|>1,$
\begin{equation}\begin{aligned}\label{e:qb1}
& |q^b_{j+1}(t,x,y)|\\
\leq& C_{14} \left(AC_{10}C_{12} \right)^j
 \,
\int_0^t\int_{\R^d} h(t-s, x,  z)  |\S^{b}_z p_0(s,z,y)|\,dz\,ds\\
\leq& C_{14} \left(AC_{10}C_{12}\right)^j
\, C_{10}A \,
\int_0^t\int_{\R^d} h(t-s,  x, z)
f_0 (s, z, y)\,dz\,ds\\
\leq& C_{14} \left(AC_{10}C_{12}\right)^{j+1}h(t, x, y).
\end{aligned}\end{equation}
On the other hand, by \eqref{e:qb1}, \eqref{e:4.5}with $\lambda=t^{-1}$  and
 $\|b^{(t^{-1})}\|_\infty=t^{1-\beta/2}\|b\|_\infty\leq t^{1-\beta/2} A$,
 noting that $|x/t^{1/2}-y/t^{1/2}|>1$ for $|x-y|>t^{1/2}$,
we have for $t^{1/2}<|x-y|\leq 1$,
$$\begin{aligned}
|q_{j+1}^b(t,x,y)|&=t^{-d/2}|q_{j+1}^{b^{(t^{-1})}} (1, x/t^{1/2}, y/t^{1/2})|\\
&\leq t^{-d/2} C_{14} (t^{1-\beta/2}AC_{10}C_{12})^{j+1} h(1, x/t^{1/2}, y/t^{1/2})\\
&\leq  C_{14} (AC_{10}C_{12})^{j+1} t^{-d/2}t^{1-\beta/2} h(1, x/t^{1/2},y/t^{1/2})\\
&\leq  C_{14} (AC_{10}C_{12})^{j+1} h(t, x, y).
\end{aligned}$$
Therefore, the above two displays prove that \eqref{k} holds for $n=j+1$ and thus for every $n\geq 1$. \qed

\begin{lem}\label{con}
For every $n\geq 0$, $q_n^b(t,x,y)$ is jointly continuous on  $(0,1]\times\R^d\times\R^d.$
\end{lem}

\pf  We prove it by induction.
Clearly $q_0^b(t,x,y)$ is continuous on $(0, 1]\times\R^d\times\R^d$.
Suppose that $q_n^b(t, x, y)$ is continuous on
$(0, 1]\times\R^d\times\R^d$.
For every $M\geq 2$ ,
it follows from \eqref{L} , Lemma \ref{L:3.3} and the dominated convergence theorem that for $\ee<1/(2M),$
$$ (t, x, y) \mapsto \int_\ee^{t-\ee}\int_{\R^d} q_n^b(t-s,x,z)\S^b_z p_0(s,z,y) dz ds
$$
is jointly continuous on $[1/M, 1]\times \R^d\times \R^d$.
On the other hand,  it follows from  \eqref{L}
that
 \begin{eqnarray*}
&&\sup_{t\in [1/M,1]} \sup_{x,y}\int_{t-\ee}^t\int_{\R^d} h(t-s,x,z)
|\S^b_z p_0(s,z,y)|\,dz\,ds\\
&\leq &c_1 A \sup_{t\in [1/M, 1]}(t-\ee)^{-(d+\beta)/2}
\sup_{x\in \R^d}\int_{t-\ee}^t\int_{\R^d} h(t-s, x, z) \,dz\,ds \\
&\leq & c_2 A (2M)^{(d+\beta)/2}  \ee ,
\end{eqnarray*}
which goes to zero as $\ee\rightarrow 0$; while
by \eqref{L} and \eqref{L1},
\begin{eqnarray*}
&& \sup_{t\in [1/M, 1]} \sup_{x,y}\int_0^\ee\int_{\R^d} h(t-s, x, z)|\S^b_z p_0 (s,z,y)|\,dz\,ds \nonumber \\
&\leq &  c_3 \left(\sup_{t\in [1/M, 1]} (t-\ee)^{-d/2}\right)
\sup_{y\in \R^d} \int_0^\ee\int_{\R^d} |\S^b_z p_0 (s,z,y)|\,dz\,ds \nonumber \\
&\leq &  c_4 (2M)^{d/2} \, \| b\|_\infty \,
\ee^{1-\beta/2}\rightarrow 0
\end{eqnarray*}
as $\ee\rightarrow 0.$
We conclude from Lemma \ref{L:3.3} and the above argument that
$$ q^b_{n+1}(t,x,y)=\int_0^t\int_{\R^d} q^b_{n}(t-s,x,z)
\S^b_z p_0(s,z,y)\,dz\,ds
$$
is jointly continuous in $(t, x, y)\in [1/M, 1]\times\R^d\times\R^d$
and so in $(t, x, y)\in (0, 1]\times\R^d\times\R^d$ by the arbitrariness of $M$.
This completes the proof of the lemma. \qed

\bigskip

\begin{lem} \label{L:3.5}
  There is a constant $C_{15} =C_{15} (d, \beta)>0$
so that for every $A>0$ and every bounded function $b$ on $\R^d\times\R^d$ with $\|b\|_\infty\leq A$ and for every integer $n\geq 0$
and $\ee >0$,
\begin{equation}\label{e:3.6b}
 \left|\int_{\{z\in \R^d: |z|>\ee\}}
\left( q^b_n (t, x+z, y) -q^b_n(t, x, y)\right)
\frac{  b(x, z)}{|z|^{d+\beta}} dz \right|
\leq  (C_{15}A)^{n+1} f_0(t, x, y)
\end{equation}
for $(t, x, y) \in (0, 1]\times \R^d\times \R^d$,
 and $\S^b_x q^b_n (t, x, y)$ exists pointwise for $(t, x, y)
\in (0, 1]\times \R^d\times \R^d$ in the sense of
\eqref{e:1.5a} with
\begin{equation}\label{e:qbn1}
\S^b_x q^b_{n+1}(t, x, y) = \int_0^t \int_{\R^d}
\S^b_x q^b_n (t-s, x, z) \S^b_z p_0(s, z, y)  dz ds
\end{equation}
and
\begin{equation}\label{e:3.11}
|\S^b_x q^b_n (t, x, y)|\leq  (C_{15}A)^{n+1} f_0(t, x, y)
\qquad \hbox{on } (0, 1]\times \R^d\times \R^d.
\end{equation}
Moreover,
\begin{equation}\label{e:qbn2}
q^b_{n+1}(t, x, y)=\int_0^t\int_{\R^d} p_0(t-s, x, z)\S^b_z q^b_{n}(s, z, y)
  dz ds
 \qquad \hbox{for } (t, x, y)\in (0, 1]\times \R^d\times \R^d.
\end{equation}
\end{lem}

\pf Let $q(t, x, y)$ denote the transition density function
  of the symmetric $\beta$-stable process on $\R^d$.
Then by \cite{CK},
we have
\begin{equation}\label{e:3.10}
q(t, x, y)\asymp t^{-d/\beta} \left( 1 \wedge \frac{t^{1/\beta}}{|x-y|} \right)^{d+\beta}\qquad \hbox{on } (0, \infty) \times \R^d \times \R^d.
\end{equation}
Observe that \eqref{e:3.9} and \eqref{e:3.10} yield
\begin{equation}
 f_0(t, x, y) \asymp  t^{-\beta/2} q(t^{\beta/2}, x, y)
 \qquad \hbox{on } (0, \infty)\times \R^d \times \R^d.
\end{equation}
Hence on $(0, \infty)\times \R^d\times \R^d$,
\begin{eqnarray*}
&& \int_0^t  \int_{\R^d} f_0(t-s, x, z)   f_0(s, z, y)  ds dz\\
&\asymp&  \int_0^t  (t-s)^{-\beta/2} s^{-\beta/2}
\left(\int_{\R^d} q((t-s)^{\beta/2}, x, z)
  q( s^{\beta/2}, z, y)  dz \right) ds \\
&=& \int_0^t  (t-s)^{-\beta/2} s^{-\beta/2}
q( (t-s)^{\beta/2}+s^{\beta/2}, x, y) ds \\
&\asymp& q( t^{\beta/2}, x, y) \int_0^t  (t-s)^{-\beta/2}
   s^{-\beta/2} ds  \\
&=& q( t^{\beta/2}, x, y)\, t^{1-(2\beta/2)} \int_0^1 (1-u)^{-\beta/2} u^{-\beta/2} du \\
&\asymp &  t^{1-\beta/2} f_0(t, x, y).
\end{eqnarray*}
In the second $\asymp$ above, we used the fact that
$$
(t/2)^{\beta/2}\leq (t-s)^{\beta/2}+s^{\beta/2}
\leq 2t^{\beta/2} \qquad \hbox{for every } s\in (0, t)
$$
and the estimate \eqref{e:3.10}, while in the last equality, we
used a change of variable $s=tu$. So there is a constant
$c_1=c_1(d, \beta)>0$ so that
\begin{equation}\label{e:3.13}
\int_0^t  \int_{\R^d} f_0(t-s, x, z)   f_0(s, z, y)  ds dz
\leq c_1 \, f_0(t, x, y) \qquad \hbox{for every } t\in (0, 1]
\hbox{ and } x, y\in \R^d.
\end{equation}
By increasing the value of $c_1$ if necessary, we may  do and assume
that $c_1$ is larger than $1$.

We now proceed by induction. Let $C_{15}:=c_1C_{10}.$
Note that
\begin{equation}\label{e:f0'}
|\S^b_xp_0(t,x,y)|\leq A|\Delta^{\beta/2}_x| p_0(t,x,y) \leq C_{10}A f_0(t, x, y).
\end{equation}
When $n=0$, \eqref{e:qbn2} holds by definition.
By Lemma \ref{1},
\eqref{e:3.6b} and \eqref{e:3.11} hold for $n=0$.
Suppose that  \eqref{e:3.6b} and \eqref{e:3.11} hold for $n=j$.
Then for every $\ee >0$,
 by the definition of $q^b_{j+1}$, Lemma \ref{L:3.3}, \eqref{e:3.13} and Fubini's theorem,
\begin{eqnarray}\label{e:3.11a}
&&   \int_{\{w\in \R^d: |\omega|>\ee\}}
\left( q^b_{j+1} (t, x+w, y) -q^b_{j+1} (t, x, y)\right)
\frac{  b(x, w)} {|w|^{d+\beta}} dw  \\
&=& \int_0^t\int_{\R^d} \left(
 \int_{\{w\in \R^d: |w|>\ee\}}
\left( q^b_j (t-s, x+w, z) -q^b_j (t-s, x, z)\right)
 \frac{  b(x, w)} {|w|^{d+\beta}} dw\right)
  \nonumber \\
&&  \hskip 0.3truein \times \, \S^b_z p_0(s, z, y)\, dz ds \nonumber
\end{eqnarray}
and so
\begin{eqnarray*}
&&  \left| \int_{\{w\in \R^d: |w|>\ee\}}
\left( q^b_{j+1} (t, x+w, y) -q^b_{j+1} (t, x, y)\right)
\frac{  b(x, w)} {|w|^{d+\beta}} dw \right| \nonumber \\
&\leq& \int_0^t\int_{\R^d} (C_{15}A)^{j+1} f_0(t-s, x, z)\,
  |\S^b_z p_0(s, z, y)| dz ds  \\
&\leq & \int_0^t\int_{\R^d} (C_{15}A)^{j+1} f_0(t-s, x, z) \, C_{10}A  f_0(s, z, y)   dz ds \\
&\leq & (C_{15}A)^{j+2} f_0(t, x, y).
 \end{eqnarray*}
By \eqref{e:3.11a} and Lebesgue dominated convergence theorem,
we conclude that
\begin{eqnarray*}
&& S^b_x q^b_{j+1}(t, x, y) \\
&:=& \lim_{\ee \to 0}
 \int_{\{w\in \R^d: |w|>\ee\}}
\left( q^b_{j+1} (t, x+w, y) -q^b_{j+1} (t, x, y)\right)
\frac{  b(x, w)} {|w|^{d+\beta}} dw \\
&=&\int_0^t\int_{\R^d} \left(
 \lim_{\ee \to 0} \int_{\{w\in \R^d: |w|>\ee\}}
\left( q^b_j (t-s, x+w, z) -q^b_j (t-s, x, z)\right)
 \frac{  b(x, w)} {|w|^{d+\beta}} dw\right)
  \\
&& \hskip 0.3truein \times \, \S^b_z p_0(s, z, y)\, dz ds \\
&=& \int_0^t\int_{\R^d} \S^b_x q^b_j (t-s, x, z)\,
     \S^b_z p_0(s, z, y)\, dz ds
\end{eqnarray*}
exists and \eqref{e:qbn1} as well as \eqref{e:3.11} holds for $n=j+1$.
(The same proof verifies \eqref{e:qbn1} when $n=0$.)
On the other hand, in view of \eqref{e:3.11} and \eqref{e:qbn2}
for $n=j$, we have by the Fubini theorem,
\begin{eqnarray*}
&&  q^b_{j+1} (t, x, y) \\
&=& \int_0^t \int_{\R^d}
q^b_j (s, x, z) \S^b_z p_0(t-s, z, y)  dz ds\\
&=& \int_0^t \int_{\R^d}
\left( \int_0^{s} \int_{\R^d} p_0(r, x, w) \S^b_w q^b_{j-1}
(s-r, w, z) dr dw \right) \S^b_z p_0(t-s, z, y)  dz ds \\
&=& \int_0^t \int_{\R^d} p_0(r, x, w)
\left( \int_r^t \int_{\R^d} \S^b_w q^b_{j-1} (s-r, w, z)
\S^b_z p_0(t-s, z, y) ds dz \right)  dw dr \\
&=& \int_0^t \int_{\R^d} p_0(r, x, w)
\S^b_w q^b_j(t-r, w, y)   dw dr . \\
\end{eqnarray*}
This verifies that \eqref{e:qbn2} also holds for $n=j+1$.
The lemma is now established by induction.
\qed

\bigskip

\begin{lem}\label{u1}
There is a  positive constant
$A_0=A_0(d,  \beta)$
so that if $\| b\|_\infty
\leq A_0$, then for every integer $n\geq 0$,
\begin{equation}\label{e:3.8}
|q_{n}^b(t,x,y)|
\leq C_{14} 2^{-n} \,  h(t, x, y)
\quad \hbox{for } t\in(0,1] \hbox{ and }  x,y\in\R^d,
\end{equation}
\eqref{e:3.6b} holds and so
$\S^b_x q_n^b(t,x,y)$ exists pointwise in the sense of \eqref{e:1.5a}
 with
\begin{equation}\label{e:3.15}
  |\S^b_x q_n^b(t,x,y)|
  \leq   2^{-n} \, f_0(t,x,y) \quad \hbox{for } t\in(0,1] \hbox{ and }  x,y\in\R^d,
\end{equation}
and
\begin{equation}\label{e:3.17}
  \sum_{n=0}^\infty q_n^b(t,x,y)\geq \frac{1}{2}p_0(t,x,y) \quad \hbox{for } t\in(0,1]
\hbox{ and } |x-y|\leq 3 t^{1/2} .
\end{equation}
\end{lem}

\pf  We take a positive constant $A_0$ so that $A_0\leq 1\wedge[2C_{10}C_{12}+2C_{15}]^{-1}.$
  We have by Lemma \ref{L:3.3} and Lemma \ref{L:3.5} that
for every $b$ with $\| b\|_\infty \leq A_0$,
$$
|q_{n}^b(t,x,y)|\leq
 C_{14}2^{-n}h(t, x, y)
\quad \hbox{ and } \quad |\S^b_x q^b_n(t, x, y)|\leq 2^{-n} f_0(t, x, y)
$$
for every $t\in (0,1]$ and  $x,y\in\R^d$.
This  establishes \eqref{e:3.8}
 and \eqref{e:3.15}.

On the other hand, by  \eqref{e:3.3}, there exists $c\geq 1$  so that $h(t, x, y)\leq cp_0(t,x,y)\asymp t^{-d/2}$
for $|x-y|\leq 3t^{1/2}$ and $t\in (0,1]$.
Take $A_0$ small enough so that $\sum_{n=1}^\infty (A_0C_{10}C_{12})^n\leq \frac{1}{2cC_{14}}.$
 Then for every $b$ with $\|b\|_\infty \leq A_0$,
 we have by Lemma \ref{L:3.3} for $|x-y|\leq 3t^{1/2}$ and $t\in (0,1]$,
$$
\sum_{n=1}^\infty |q^b_n(t,x,y)|\leq    cC_{14}\sum_{n=1}^\infty (A_0C_{10}C_{12})^n
 p_0(t,x,y)\leq \frac{1}{2} p_0(t,x,y).
$$
Consequently, for $|x-y|\leq 3t^{1/2}$ and $t\in (0,1]$,
$$
\sum_{n=0}^\infty q_n^b(t,x,y)\geq p_0(t,x,y)-\sum_{n=1}^\infty |q_n^b(t,x,y)|\geq \frac{1}{2}p_0(t,x,y).
$$
\qed

\medskip

We now extend the results in Lemma \ref{u1}
to any bounded $b$ that satisfies condition \eqref{e:b}.
Recall that $A_0$ is the positive constant in Lemma \ref{u1}.

\medskip

\begin{thm}\label{T:3.7}
For every  $A>0$, there are positive constants $C_{16}=C_{16}(d, \beta, A)>1$ and $0<C_{17}=C_{17}(d, \beta, A)<1$
 so that for every bounded function $b$ with $\| b\|_\infty \leq A$,
  that satisfies condition
\eqref{e:b} and $n\geq 0$,
\begin{equation}\label{e:3.16n}
  |q_n^b(t,x,y)|
  \leq C_{16} 2^{-n} p_1(t,C_{17}x,C_{17}y)
\end{equation}
for every $0<t \leq 1\wedge ( A_0/\|b\|_\infty)^{2/(2 -\beta)}$
and $x,y\in\R^d$, and
\begin{equation}\label{e:3.17n}
  \sum_{n=0}^\infty q_n^b(t,x,y)\geq \frac{1}{2}p_0(t,x,y) \quad \hbox{for }
  0<t \leq 1\wedge( A_0/\|b\|_\infty)^{2/(2 -\beta)}
\hbox{ and } |x-y|\leq 3 t^{1/2} .
\end{equation}
Moreover, for every $n\geq 0$, \eqref{e:3.6b} holds and so
$\S^b_x q_n^b(t,x,y)$ exists
pointwise in the sense of \eqref{e:1.5a} with
\begin{equation}\label{e:3.24}
  |\S^b_x q_n^b(t,x,y)|
  \leq  2^{-n}  f_0(t, x, y)
\end{equation}
for every $0<t \leq 1\wedge ( A_0/\|b\|_\infty)^{2/(2 -\beta)}$
and $x,y\in\R^d$. Moreover, \eqref{e:qbn1} and \eqref{e:qbn2} hold.
\end{thm}

\pf Note that there exist $c_k>0, k=1,2$ such that $h(t,x,y)\leq c_1 p_1(t, c_2x, c_2y)$ for $t\in (0,1]$ and $x,y\in\R^d.$
Thus in view of Lemma \ref{u1}, it suffices to prove the theorem
 for $A_0<\| b\|_\infty \leq A$.
Set $r=( \|b\|_\infty/A_0)^{2/(2 -\beta)}$.
The function $b^{(r)}$ defined by \eqref{e:4.1} has
the property $\|b^{(r)} \|_\infty =A_0$.
Thus by Lemma \ref{u1},
  for every integer $n\geq 0$,
\begin{equation}\label{e:3.8c}
|q_{n}^{b^{(r)}}(t,x,y)|
\leq C_{14} \, 2^{-n} \,  h(t, x, y)
\quad \hbox{for } t\in(0,1] \hbox{ and }  x,y\in\R^d,
\end{equation}
\eqref{e:3.6b} holds and so
$\S^b_x q_n^{b^{(r)}}(t,x,y)$ exists pointwise in the sense of \eqref{e:1.5a}
 with
\begin{equation}\label{e:3.15c}
  |\S^b_x q_n^{b^{(r)}}(t,x,y)|
  \leq   2^{-n} \, f_0(t,x,y) \quad \hbox{for } t\in(0,1] \hbox{ and }  x,y\in\R^d,
\end{equation}
and
\begin{equation}\label{e:3.17c}
  \sum_{n=0}^\infty q_n^{b^{(r)}}(t,x,y)\geq \frac{1}{2}p_0(t,x,y) \quad \hbox{for } t\in(0,1]
\hbox{ and } |x-y|\leq 3 t^{1/2} .
\end{equation}
We have
by   \eqref{e:4.5}, \eqref{e:3.8c} and \eqref{e:1.0} that
for every $0<t\leq 1/r= (A_0/\|b\|_\infty)^{2/(2 -\beta)}$
and $x, y\in \R^d$,
\begin{eqnarray*}
|q_n^b(t,x,y)|&=&
r^{d/2} \,  |q^{b^{(r)}}_n  (r t, r^{1/2} x, r^{1/2} y) |  \\
&\leq& C_{14} 2^{-n} r^{d/2} \, h
 (r t, r^{1/2} x, r^{1/2} y)  \\
 &\leq & C_{14}\, 2^{-n} \left( t^{-d/2}\wedge
  \left( p_0(t, x, y) + \frac{ r^{1-\beta/2}t}{|x-y|^{d+\beta}} \right)\right)\\
&\leq& C_{16}\, 2^{-n} p_1(t, C_{17}x, C_{17}y),
\end{eqnarray*}
which establishes \eqref{e:3.16n}.
Similarly,  \eqref{e:3.17n}  follows
from \eqref{e:4.3}, and \eqref{e:3.17c}, while
the conclusion of \eqref{e:3.24} is a direct consequence of
\eqref{e:4.2}, \eqref{e:4.5} and \eqref{e:3.15c}.
That \eqref{e:qbn1} and \eqref{e:qbn2} hold follows directly
from Lemma \ref{L:3.5} and Lemma \ref{L:3.5b}.
\qed

\bigskip

Recall that  $q^b(t,x,y):=\sum_{n=0}^\infty q_n^b(t,x,y)$, whenever
it is convergent.
The following theorem follows immediately from Lemmas \ref{con},
 \ref{u1} and Theorem\ref{T:3.7}.

\begin{thm}\label{c1}
For every $A>0$,
let $C_{16}=C_{16}(d, \beta, A)$ and $C_{17}=C_{17}(d, \beta, A)$
 be the constants in Theorem
\ref{T:3.7}. Then for every bounded function
$b$ with $\| b\|_\infty \leq A$ that satisfies condition
\eqref{e:b}, $q^b(t,x,y)$ is well defined and is
jointly continuous in $(0, 1\wedge ( A_0/\|b\|_\infty)^{2/(2 -\beta)} ]\times\R^d\times\R^d$.
Moreover,
$$
  |q^b(t, x, y)| \leq 2C_{16}\, p_1(t,C_{17}x,C_{17}y)
$$
and $\S^b_x q^b(t, x, y)$ exists pointwise in the sense of
\eqref{e:1.5a} with
$$ | \S^b_x q^b(t, x, y) | \leq  2 f_0(t, x, y)
$$
  for every $0<t \leq 1\wedge ( A_0/\|b\|_\infty)^{2/(2 -\beta)}$
and $x,y\in\R^d$, and
\begin{equation}\label{e:ql1}
q^b(t,x,y)\geq \frac{1}{2}p_0(t,x,y) \quad \hbox{for }
  0<t \leq 1\wedge ( A_0/\|b\|_\infty)^{2/(2 -\beta)}
\hbox{ and } |x-y|\leq 3 t^{1/2} .
\end{equation}
Moreover, for every
$0<t \leq 1\wedge ( A_0/\|b\|_\infty)^{2/(2 -\beta)}$ and $x,y\in\R^d$,
\begin{eqnarray}\label{e:3.21}
q^b(t, x, y)&=&p_0(t, x, y)+\int_0^t \int_{\R^d} q^b(t-s, x, z)
\S^b_z p_0(s, z, y)  dz ds\\
&=& p_0(t, x, y)+\int_0^t \int_{\R^d} p_0 (t-s, x, z)
\S^b_z q^b (s, z, y)dz ds \label{e:Du2}.
\end{eqnarray}
\end{thm}

\bigskip

\begin{thm}\label{T:3.9}
Suppose that $b$ is a bounded function on $\R^d\times \R^d$ satisfying \eqref{e:b}. Let $A_0$ be the constant in Lemma \ref{u1}.
 Then for every  $t, s>0$
with $t+s\leq 1\wedge ( A_0/\|b\|_\infty)^{2/(2 -\beta)}$ and $x, y\in \R^d$,
\begin{equation}\label{e:3.26}
\int_{\R^d} q^b(t, x, z) q^b(s, z, y) dz = q^b(t+s, x, y).
\end{equation}
\end{thm}

\pf In view of Theorem \ref{T:3.7}, we have
$$ \int_{\R^d} q^b(t, x, z) q^b(s, z, y) dz
=\sum_{j=0}^\infty \sum_{n=0}^j \int_{\R^d} q^b_n(t, x, z)
q^b_{j-n}(s, z, y) dz.
$$
So it suffices to show that for every $j\geq 0$,
\begin{equation}\label{e:3.12}
\sum_{n=0}^j \int_{\R^d} q^b_n(t, x, z)
q^b_{j-n}(s, z, y) dz=q^b_j (t+s, x, y)
\end{equation}
Clearly, \eqref{e:3.12} holds for $j=0$.
Suppose that \eqref{e:3.12} holds for $j=l\geq 1$.
Then we have by Fubini's theorem and the estimates
in \eqref{L} and Theorem \ref{T:3.7},
\begin{eqnarray*}
&& \sum_{n=0}^{l+1} \int_{\R^d} q^b_n(t, x, z)
   q^b_{l+1-n}(s, z, y) dz \\
&=&  \int_{\R^d} q^b_{l+1}(t, x, z)
   p_0 (s, z, y) dz + \sum_{n=0}^{l} \int_{\R^d} q^b_n(t, x, z)
   q^b_{l+1-n}(s, z, y) dz  \\
&=& \int_{\R^d} \left(   \int_0^t \int_{\R^d}
     q^b_{l}(t-r, x, w) \S^b_w p_0(r, w, z)  dw dr \right)
     p_0 (s, z, y) dz \\
&& + \sum_{n=0}^{l} \int_{\R^d} q^b_n(t, x, z)
   \left( \int_0^s\int_{\R^d} q^b_{l-n}(s-r, z, w)\S^b_w
   p_0(r, w, y)  dw dr\right)  dz \\
&=& \int_0^t \int_{\R^d}  q^b_{l}(t-r, x, w) \S^b_w
    p_0(r+s, w, y)  dw dr\\
&& +  \int_0^s  \int_{\R^d} q^b_l(t+s-r, x, w) \S^b_w
   p_0(r, w, y) dw dr   \\
&=& q^b_{l+1} (t+s, x, y).
\end{eqnarray*}
This proves that \eqref{e:3.12} holds for $j=l+1$.
So by induction, we conclude that \eqref{e:3.12} holds for
every $j\geq 0$.  \qed

\bigskip

For notational simplicity,
denote $1\wedge ( A_0/\|b\|_\infty)^{2/(2 -\beta)}$
 by $\delta_0$.
 In view of Theorem \ref{T:3.9}, we can uniquely extend the definition of
$q^b(t, x, y)$ to $t> \delta_0$
by using the Chapman-Kolmogorov equation recursively as follows.

Suppose that $q^b(t, x, y)$ has been defined and satisfies the
Chapman-Kolmogorov equation \eqref{e:3.26} on $(0, k\delta_0]\times \R^d
\times \R^d$. Then for $t\in (k\delta_0, (k+1)\delta_0]$, define
\begin{equation}\label{pk}
q^b(t,x,y) =\int_{\R^d}q^b(s,x,z)q^b(r,z,y)\,dz, \quad  x, y\in \R^d
\end{equation}
for any $s, r\in (0, k\delta_0]$ so that $s+r=t$.
Such $q^b(t, x, y)$ is well defined on $(0, \infty)\times \R^d\times \R^d$
and satisfies \eqref{e:3.26} for every $s, t>0$.
Moreover, since Chapman-Kolmogorov equation holds for $q^b(t, x, y)$ for all $t, s>0$,
we have by Theorem \ref{c1}
that for every $A>0$,
there are constants
$c_i=c_i(d, \beta, A)$, $i=1, 2, 3$, so that for every
$b(x, z)$ satisfying \eqref{e:b} with $\| b\|_\infty \leq A$,
\begin{equation}\label{e:3.26n}
|q^b(t, x, y)| \leq c_1 \, e^{c_2t} \,   p_1 (t, c_3x, c_3y)
\qquad \hbox{for every } t>0 \hbox{ and } x, y\in \R^d .
\end{equation}

By the induction method and the Chapman-Kolmogrov equation for $q^b(t,x,y),$  \eqref{e:3.21} and \eqref{e:Du2} can be extended
to every $t>0$ and $x,y\in\R^d.$

\begin{thm}\label{T:3.10a}
$q^b(t, x, y)$ satisfies \eqref{e:3.21} and \eqref{e:Du2} for every $t>0$ and $x,y\in\R^d.$
\end{thm}

\pf Let $\delta_0:=1\wedge ( A_0/\|b\|_\infty)^{2/(2 -\beta)}.$
It suffices to prove that for every $n\geq 1,$
\eqref{e:3.21} and \eqref{e:Du2}  hold for all $t\in (0,n\delta_0]$ and $x,y\in\R^d.$

Clearly, \eqref{e:3.21} holds for $t\in (0,n\delta_0]$ with $n=1.$
Suppose that  \eqref{e:3.21} holds for $t\in (0,n\delta_0]$ with $n=k.$
For $t\in (k\delta_0, (k+1)\delta_0]$, take $l,s\in (0, k\delta_0]$ so that $l+s=t.$
Then we have by Fubini's theorem, Chapman-Kolmogorov equation of $q^b$, Lemma \ref{2'}, \eqref{L} and (\ref{e:3.26n}),
 \begin{eqnarray*}
 q^b(l+s, x,y)
&=& \int_{\R^d} q^b(l,x,z)q^b(s,z,y)\,dz\\
&=& \int_{\R^d} q^b(l,x,z) \left(p_0(s,z,y)+\int_0^s\int_{\R^d} q^b(s-r, z, \omega) \S^b_\omega p_0(r, \omega, y)\,d\omega\,dr\right)\,dz\\
&=& \int_{\R^d} p_0(l,x,z) p_0(s,z,y)\,dz\\
&& +\int_{\R^d} \left(\int_0^l\int_{\R^d} q^b(l-u, x, \eta) \S^b_\omega p_0(u, \eta, z)\,d\eta\,du\right) p_0(s,z,y)\,dz\\
&& +\int_0^s\int_{\R^d} q^b(l+s-r, x, \omega) \S^b_\omega p_0(r, \omega, y)\,d\omega\,dr\\
&=& p_0(l+s,x,y)+\int_0^l\int_{\R^d} q^b(l-u, x, \eta) \S^b_\omega p_0(u+s, \eta, y)\,d\eta\,du\\
&& +\int_0^s\int_{\R^d} q^b(l+s-r, x, \omega) \S^b_\omega p_0(r, \omega, y)\,d\omega\,dr\\
&=& p_0(l+s,x,y)+\int_0^{l+s} \int_{\R^d} q^b(l+s-r, x, z)\S^b_z p_0(r,z,y)\,dz\,dr.
\end{eqnarray*}
By the similar procedure as above, we can also prove that \eqref{e:Du2} holds for every $t>0$ and $x,y\in\R^d.$
\qed

\bigskip

\begin{thm}\label{T:3.10}
Suppose that $b$ is a bounded function on $\R^d\times \R^d$ satisfying \eqref{e:b}.
$q^b(t, x, y)$ is the unique continuous kernel that satisfies
the Chapman-Kolmogorov equation \eqref{e:3.26} on $(0, \infty)\times \R^d
\times \R^d$ and that
 for some $\ee>0, c_1>1$ and $0<c_2<1$,
\begin{equation}\label{e:3.22}
|q^b(t, x, y)| \leq c_1 \, p_1 (t, c_2x, c_2y)
\end{equation}
 and \eqref{e:3.21} hold
 for $(t, x, y)\in (0, \ee]\times \R^d \times \R^d$.
 Moreover, \eqref{e:3.26n} holds for $q^b(t, x, y)$.
\end{thm}

\pf Suppose that $\overline  q$ is any continuous kernel that satisfies,
for some $\ee>0$,  \eqref{e:3.21} and \eqref{e:3.22} hold for
$(t, x, y)\in (0, \ee]\times \R^d \times \R^d$.
Without loss of generality, we may and do assume that
$\ee <1\wedge (A_0/\| b\|_\infty)^{2/(2 -\beta)}$.
Using \eqref{e:3.21} recursively, one gets
\begin{equation} \label{e:3.29}
\overline  q (t, x, y)=\sum_{j=0}^n q^b_j (t, x, y)
 + \int_0^t \int_{\R^d}  \overline  q (t-s, x, z)
 (\S^b p_0)^{*, n+1}_z (s,z,y)
ds dz.
 \end{equation}
Here $(\S^b p_0)^{*, n}_z (s,z,y) $  denotes the $n$th
convolution operation of the function $\S^b_z p_0 (s, z, y)$;
that is, $(\S^b p_0)^{*, 1}_z (s,z,y)= \S^b_z p_0 (s,z,y)$
and
\begin{equation}\label{e:s}
(\S^b p_0)^{*, n}_z (s,z,y)= \int_0^s \int_{\R^d}
\S^b_z p_0(r,z,w) \, (\S^b p_0)^{*, n-1}_w (s-r,w,y) dw dr
\quad \hbox{for } n\geq 2.
\end{equation}
It follows from \eqref{e:f0'} and \eqref{e:3.13} that for every $A>0$ so that $\|b\|_\infty\leq A,$
$$
|(\S^b p_0)^{*, n}_z (s,z,y)| \leq (C_{15}A)^{n} f_0(t, x, y),
$$
where $C_{15}$ is the constant in Lemma \ref{L:3.5}.
Noting that the constant $A_0$ defined in Lemma \ref{u1} satisfies  $A_0\leq 1/(2C_{15}).$
So for every bounded function $b$ with $\|b\|_\infty\leq A_0,$ we have
\begin{equation}\label{e:f1}
|(\S^b p_0)^{*, n}_z (s,z,y)|\leq 2^{-n} f_0(s,z,y), \quad s\in (0, 1).
\end{equation}
By the scale change formulas \eqref{e:4.2} and \eqref{e:4.3}, when $\|b\|_\infty>A_0,$
$$
|(\S^b p_0)^{*, n}_z (s,z,y)|\leq 2^{-n} f_0(s,z,y), \quad s\in (0, 1\wedge (A_0/\| b\|_\infty)^{2/(2 -\beta)}).
$$
Thus, by the condition
\eqref{e:3.22} and Lemma \ref{2'},
$$ \left| \int_0^t \int_{\R^d}  \overline q (t-s, x, z)
(\S^b p_0)^{*, n}_z (s,z,y) ds dz \right| \leq c_3 2^{-n} p_1(t,c_4x,c_4y).
$$

It follows from \eqref{e:3.29} that
$$  \overline q (t, x, y)=\sum_{n=0}^\infty q_n^b (t,x, y)
=q^b(t, x, y)
$$
for every $t\in (0, \ee]$ and $x, y\in \R^d$.
Since both $ \overline q$ and $q^b$ satisfy the Chapman-Kolmogorov
equation \eqref{e:3.26},
$\overline q=q^b$ on $(0, \infty) \times \R^d \times \R^d$.
\qed

\bigskip

\begin{remark}\rm
It follows from the definition of $q^b_n(t, x, y)$ and
Lemma \ref{L:3.5} that
$ (\S^b p_0)^{*, n+1} (s,z,y) = \S^b_z q^b_n(s, z, y)$.
\qed
\end{remark}

\bigskip

\bigskip

In view of Lemma \ref{L:3.5b} and Chapman-Kolmogorov equation, we have

\begin{thm}\label{T:3.11}
Suppose that $b$ is a bounded function on $\R^d\times \R^d$ satisfying \eqref{e:b}.
 \ $q^b(t, x, y)=\lambda^{d/2} q^{b^{(\lambda)}}(\lambda t, \lambda^{1/2} x, \lambda^{1/2} y)$ on $(0, \infty)\times
\R^d\times \R^d$, where $b^{(\lambda)}(x, z):=\lambda^{\beta/2-1}
 b(\lambda^{-1/2} x, \lambda^{-1/2} z)$.
\end{thm}

\bigskip

For a bounded function $f$ on $\R^d$, $t>0$ and $x\in\R^d$,
 we define
$$
T^b_tf(x)=\int_{\R^d} q^b(t,x,y)f(y)\,dy
\quad \hbox{and} \quad
P_t f(x)=\int_{\R^d} p_0(t, x, y) f(y) dy.
$$
The following lemma follows immediately from \eqref{e:3.26}
and \eqref{pk}.

\medskip

\begin{lem}\label{lq^b}
Suppose that $b$ is a bounded function on $\R^d\times \R^d$ satisfying \eqref{e:b}.
For all $s,t>0,$ we have $T^b_{t+s}=T^b_tT^b_s.$
\end{lem}

\medskip

\begin{thm}\label{T:3.13}
Let $b$ be a bounded function on $\R^d\times \R^d$ satisfying \eqref{e:b}.
Then for every $f\in C^2_b (\R^d)$,
$$
T^b_t f(x) -f(x)=\int_0^t T^b_s \L^b f(x) ds
\qquad \hbox{for every } t>0, \, x\in \R^d .
$$
\end{thm}

\pf  Note that by Theorem \ref{T:3.10a}, for $f\in C^2_b (\R^d)$,
\begin{equation}\label{e:3.32}
T^b_t f(x)=P_t f(x)+\int_0^t T^b_{t-s} \S^b P_s f(x) ds
=P_t f(x)+\int_0^t T^b_s \S^b P_{t-s} f(x) ds.
\end{equation}
Hence
\begin{eqnarray*}
&& T^b_t f(x)-f(x) \\
&=&P_t f(x)-f(x) +\int_0^t T^b_s \S^b f(x) ds +
\int_0^t T^b_s \S^b (P_{t-s} f-f) (x) ds \\
&=& \int_0^t P_s \Delta f(x) ds
+ \int_0^t T^b_s \S^b f(x) ds +
\int_0^t T^b_s \S^b (P_{t-s} f-f) (x) ds      \\
&=& \int_0^t T^b_s \Delta f(x) ds
 - \int_0^t \left(\int_0^{s} T^b_r \S^b P_{s-r}(\Delta f)(x) dr \right) ds \\
&& + \int_0^t T^b_s \S^b f(x) ds +
\int_0^t T^b_s \S^b (P_{t-s} f-f) (x) ds  \\
&=& \int_0^t T^b_s \left( \Delta+\S^b\right) f(x) ds
 - \int_0^t \left(\int_r^{t} T^b_r \S^b P_{s-r}(\Delta f)(x) ds
  \right) dr \\
&& +
\int_0^t T^b_s \S^b (P_{t-s} f-f) (x) ds  \\
&=& \int_0^t T^b_s  \L^b  f(x) ds
 - \int_0^t   T^b_r \S^b ( P_{t-r}f -f)(x)   dr  +
\int_0^t T^b_s \S^b (P_{t-s} f-f) (x) ds  \\
&=& \int_0^t T^b_s  \L^b  f(x) ds.
\end{eqnarray*}
Here in the third inequality, we used \eqref{e:3.32};
while in the fifth inequality we used Lemma \ref{1}
and \eqref{e:3.26n}, which
allow the interchange of the integral sign $\int_r^t$
with $T^b_r \S^b$, and
the fact that
$$  \int_r^t  P_{s-r}(\Delta f )(x)ds
= \int_r^t \left(\frac{d}{ds} P_{s-r}f (x)\right) ds
= P_{t-r}f(x)-f(x).
$$
\qed

\bigskip

\begin{thm}\label{lambda1}
Let $b$ be a bounded function on $\R^d\times \R^d$ satisfying \eqref{e:b}.
 Then $q^b(t,x,y)$ is jointly continuous in $(0,\infty)\times\R^d\times\R^d$  and $\int_{\R^d}q^b(t,x,y)\,dy=1$ for
 every $x\in \R^d$ and  $t>0.$
\end{thm}

\pf  By Lemma \ref{lq^b}, we have
\begin{equation}\label{eq^b}
q^b(t+s,x,y)=\int_{\R^d}q^b(t,x,z)q^b(s,z,y)\,dz, \quad x,y\in\R^d, s,t>0.
\end{equation}
Continuity of $q^b(t,x,y)$ in $(t, x, y)\in (0, \infty) \times \R^d
\times \R^d$ follows from Theorem \ref{c1}, (\ref{eq^b}) and the dominated convergence theorem.
For $n\geq 1$ and $t\in (0, T]$,
it follows from \eqref{L}, Lemma \ref{2'}, Theorem \ref{T:3.7} and  Fubini's Theorem
that for every $t\in (0, 1\wedge( A_0/\|b\|_\infty)^{2/(2 -\beta)}]$,
$$\begin{aligned}
\int_{\R^d}q_n^b(t,x,y)\,dy
&=\int_{\R^d}\int_{\R^d}\int_0^t q^b_{n-1}(t-s,x,z) \S_z^b p_0(s,z,y)\,ds\,dz\,dy\\
&=\int_{\R^d}\int_0^t q^b_{n-1}(t-s,x,z)\S^b_z \left(\int_{\R^d} p_0(s,z,y)\,dy\right)\,ds\,dz =0.\\
\end{aligned}$$
Hence we have by Lemma \ref{u1},
$$\int_{\R^d}q^b(t,x,y)\,dy=\int_{\R^d} p_0(t,x,y)\,dy=1$$ for
$t\in (0,1\wedge( A_0/\|b\|_\infty)^{2/(2 -\beta)}]$. This conservativeness property extends to all $t>0$ by (\ref{eq^b}). \qed

\bigskip

Theorem \ref{T:1.1} now follows from
\eqref{e:1.0}, \eqref{e:h},
Theorems \ref{c1}, \ref{T:3.10a}, \ref{T:3.10},
 \ref{T:3.13} and \ref{lambda1}.

\section{$C_\infty $-Semigroups   and Positivity }\label{S:4}

Recall that $A_0$ is the positive constant in Lemma \ref{u1}.

\begin{lem}\label{lambda2}
Suppose that $b$ is a bounded function on $\R^d\times \R^d$
satisfying condition \eqref{e:b}.
Then $\{T^b_t, t>0\}$ is a strongly continuous
  semigroup in $C_\infty(\R^d).$
\end{lem}

\pf  Note that $q^b(t,x,y)$ is jointly continuous in $(0, \infty) \times \R^d
\times \R^d$ and
 there are constants $c_k, k=1,2,3$ so that
$$
|q^b(t,x,y)|\leq c_1e^{c_2t}   p_1  (t, c_3x, c_3y)
\qquad \hbox{for every } t>0 \hbox{ and } x, y \in \R^d.
$$
The  proof is a minor modification of that for \cite[Proposition 2.3]{CKS}.
\qed

\bigskip

\begin{lem}\label{L:4.2}
Let  $b$ be a bounded function on $\R^d\times \R^d$ satisfying
\eqref{e:5.1}.
For each $f\in C^2_\infty (\R^d)$, $\L^b f(x)$ exists
pointwise and is in $C_\infty (\R^d)$.
\end{lem}

\pf Suppose that $f\in C^2_\infty (\R^d)$.
Denote $\sum_{i, j=1}^d | \partial^2_{ij} f(x)|$ by
$|D^2 f(x)|$.
Let $R>1$ to be chosen later.
Then for each $x\in \R^d$, we have by Taylor expansion,
\begin{eqnarray*}
\Phi_f (x)&:=& \int_{\R^d} \left| f(x+z)-f(x)-\nabla f (x)
\cdot z \1_{\{|z|\leq 1\}} \right|
\frac{1}{|z|^{d+\beta}} dz \non
&\leq & \int_{|z|\leq 1} \left| f(x+z)-f(x)-\nabla f (x) \cdot z
\1_{\{|z|\leq 1\}} \right|
\frac{1}{|z|^{d+\beta}} dz \non
&& + \int_{1<|z|\leq R} \left| f(x+z)-f(x) \right|
\frac{1}{|z|^{d+\beta}} dz  + \int_{|z|>R} \left| f(x+z)-f(x) \right|
\frac{1}{|z|^{d+\beta}} dz \non
&\leq &   \sup_{|y|\leq 1} | D^2 f(x+y)| \int_{|z|\leq 1} |z|^{2-d-\beta}dz
+ \int_{1<|z|\leq R} \left| f(x+z)-f(x) \right|
\frac{1}{|z|^{d+\beta}} dz \non
&&+ 2 \|f \|_\infty \int_{|z|>R} |z|^{-d-\beta} dz \non
&=& c \sup_{|y|\leq 1} | D^2 f(x+y)|+ \int_{1<|z|\leq R} \left| f(x+z)-f(x) \right|
\frac{1}{|z|^{d+\beta}} dz
 + c R^{-\beta}  \|f \|_\infty .
\end{eqnarray*}
For any given $\ee >0$, we can take $R$ large so that
$c R^{-\beta}  \|f \|_\infty <\ee/2$ to conclude
that
\begin{equation}\label{e:4.3b}
\lim_{|x|\to \infty}   \int_{\R^d} \left| f(x+z)-f(x)-\nabla f (x) \cdot z \1_{\{|z|\leq 1\}} \right|
\frac{1}{|z|^{d+\beta}} dz =0.
\end{equation}
By the same reason, applying the above argument to function
$x\mapsto f(x+y)-f(x)$ in place of $f$ yields that for every
$\ee>0$ and $x_0\in \R^d$, there is $\delta>0$ so that
\begin{equation}\label{e:4.4b}
\Phi_{f(\cdot +z)-f}(x_0)<\ee \quad \hbox{for every } |z|<\delta.
\end{equation}
It follows from the last two displays, the definition of $\L^b$ and
\eqref{e:5.1}
that $\L^bf(x)$ exists for every $x\in \R^d$ and $\L^bf\in C_\infty (\R^d)$.
\qed

\bigskip

\noindent{\bf Proof of Theorem \ref{T:1.2}.}
 Since $b$ satisfies condition \eqref{e:5.1}, it is easy to verify
that $\L^b f\in C_\infty (\R^d)$ for every $f\in C^2_c(\R^d)$.
Let $\wh \L^b$ denote the infinitesimal generator of the strongly continuous  semigroup $\{T^b_t; t\geq 0\}$ in $C_\infty (\R^d)$, which is
a closed linear operator.
It follows from Theorem \ref{T:3.13},
 Lemmas \ref{lambda2} and \ref{L:4.2}
that for every
$f\in C^2_\infty (\R^d)$,
$(T^b_t f(x)-f(x))/t$  converges uniformly to $\L^b f(x)$ as $t\to 0$.
So
\begin{equation}\label{e:4.5b}
C^2_\infty(\R^d)\subset
D(\wh \L^b) \quad \hbox{ and } \quad
\wh \L^b f=\L^b f \quad \hbox{for } f\in C^2_\infty (\R^d).
\end{equation}
In view of Theorem \ref{c1},
there are constants $c_k>0, k=1,2,3$ so that
\eqref{e:3.26n} holds. This implies that
$$
\sup_{x\in \R^d}\int_0^\infty e^{-\lambda t} |T^b_t f|(x) dt \leq c_\lambda \|f \|_\infty,
\quad f\in C_\infty (\R^d),
$$
for every $ \lambda>c_2$. Observe that $e^{-c_2 t}T^b_t$ is a strongly continuous
semigroup in $C_\infty (\R^d)$ whose infinitesimal generator is $\wh \L^b-c_2$.
The above display
implies that $(0, \infty)$ is contained in the residual set $\rho (\wh \L^b -c_2)$ of
$\wh \L^b-c_2$.  Therefore by Theorem \ref{lambda1} and the Hille-Yosida-Ray theorem \cite[p165]{EK}, $\{e^{-c_2 t}T^b_t; t\geq 0\}$
is a positive preserving semigroup on $C_\infty (\R^d)$
if and only if $\wh \L^b-c_2$ satisfies the positive maximum principle.
On the other hand,   Courr\'ege's first theorem (see \cite[p158]{D})
tells us that $\wh \L^b-c_2$ satisfies the positive maximum principle
if and only if for each $x\in \R^d$,
$$  b(x,z) \geq 0
\quad \hbox{for a.e. } z\in \R^d.
$$
Since $ e^{-c_2 t}T^b_t$ has a continuous integral kernel $e^{-c_2t}q^b(t, x, y)$, it follows that
$q^b(t, x, y)\geq 0$ on $(0, \infty)\times \R^d \times \R^d$
if and only if for each $x\in \R^d$,  \eqref{e:1.15} holds.
\qed

\section{Feller process and heat kernel estimates}\label{S:5}

Suppose that $b$ is a bounded  function
satisfying conditions \eqref{e:b}, \eqref{e:5.1} and \eqref{e:1.15}.
Then it follows from Theorem \ref{T:1.2} and
Lemma \ref{lambda2}, $T^b$ is a Feller semigroup.
So it uniquely determines a conservative
Feller process $X^b=\{X^b_t, t\geq 0, \P_x, x\in \R^d\}$ having
$q^b(t, x, y)$ as its transition density function.
Since, by Theorem \ref{T:3.10},
$q^b(t, x, y)$ is continuous and $q^b(t, x, y)\leq c_1 e^{c_2t} p_1  (t, c_3x, c_3y)$
for some positive constants $c_k, k=1,2,3$,
$X^b$ enjoys the strong Feller property.

\bigskip

\begin{prp}\label{p2}
Suppose that $b$ is a bounded  function
satisfying conditions \eqref{e:b}, \eqref{e:5.1} and \eqref{e:1.15}.
For each $x\in\R^d$ and $f\in C_b^2 (\R^d),$
$$M_t^f:=f(X^b_t)-f(X^b_0)-\int_0^t \L^b f(X^b_s)\,ds$$
is a martingale under $\P_x$.  So in particular,
the Feller process $(X^b, \P_x, x\in \R^d)$ solves the martingale
problem for $(\L^b, C^2_\infty (\R^d))$.
\end{prp}

\pf This follows immediately from
Theorem \ref{T:3.13}
and the Markov property of $X^b$. \qed

We next determine the L\'evy system of $X^b$. Recall that
\begin{equation}\label{e:5.2}
J^b(x, y)=\frac{  \, b(x, y-x)}{|x-y|^{d+\beta}}.
\end{equation}

By Proposition \ref{p2} and the similar argument in \cite[Theorem 2.6]{CKS}, we have the following result.

\begin{prp}\label{p3}
Suppose that $b$ is a bounded  function
satisfying conditions \eqref{e:b}, \eqref{e:5.1} and \eqref{e:1.15}.
Assume that $A$ and $B$ are disjoint compact sets in $\R^d$.
Then
$$
\sum_{s\leq t} \1_{\{X^b_{s-}\in A, X^b_s\in B\}}-
\int_0^t \1_A(X^b_s) \int_B J^b(X^b_s, y) dy\,ds
$$
is a $\P_x$-martingale for each $x\in \R^d$.
\end{prp}

Proposition \ref{p3} implies that
$$
\E_x\left[ \sum_{s\le t}{\bf 1}_A(X^b_{s-}){\bf 1}_B(X^b_s)
\right]=
\E_x\left[\int^t_0\int_{\R^d} {\bf 1}_A(X^b_s){\bf 1}_B(y)J^b(X^b_s, y)dyds\right].
$$
Using this and a routine measure theoretic argument, we get
$$\E_x\left[ \sum_{s\le
t}f(s, X^b_{s-}, X^b_s) \right]
=\E_x\left[\int^t_0\int_{\R^d}f(s, X^b_s, y)J^b(X^b_s, y)dyds\right]
$$
for any non-negative measurable function $f$ on $(0, \infty)
\times \R^d\times \R^d$
vanishing on $\{(x, y)\in \R^d\times \R^d: x=y\}$. Finally,
following the same arguments as in \cite[Lemma 4.7]{CK} and
\cite[Appendix A]{CK2}, we get

\begin{prp}\label{p4}
Suppose that $b$ is a bounded  function
satisfying conditions \eqref{e:b}, \eqref{e:5.1} and \eqref{e:1.15}.
Let $f$ be a nonnegative function on $\R_+\times
\R^d\times\R^d$ vanishing on the diagonal. Then for stopping time $T$
with respect to the minimal admissible filtration generated by $X^b$,
$$
\E_x \left[ \sum_{s\leq T}f(s, X^b_{s-}, X^b_s) \right]
=\E_x \left[ \int_0^T\int_{\R^d}
f(s,X^b_s,u) J^b(X^b_s,u) \,du\,ds \right].
$$
\end{prp}

 \bigskip

 To remove the assumption \eqref{e:5.1} on $b$, we approximate
 a general measurable function $b(x, z)$ by continuous $k_n(x, z)$.
 To show that $q^{k_n}(t, x, y)$ converges to $q^b (t, x, y)$,
 we establish equi-continuity of $q^b(t, x, y)$ and apply
 the uniqueness result, Theorem \ref{T:3.10}.

 \bigskip

\begin{prp}\label{H1}
 For each  $0<t_0<T<\infty$ and $A>0$,
 the function $q^b(t,x,y)$ is uniform continuous in $(t,x)\in (t_0,T)\times\R^d$ for
every $b$ with $\|b\|_\infty\leq A$ that satisfies \eqref{e:b} and for
all $y\in \R^d.$
\end{prp}

\pf In view of Theorem \ref{T:3.11},
it suffices to prove the theorem
for $A=A_0$, where $A_0$ is the constant in Lemma \ref{u1}
(or in Theorem \ref{T:1.1}).
Using the Chapman-Kolmogorov equation for $q^b(t, x, y)$ , it suffices to prove the Proposition for $T=1$.

Noting that $q^b_n, n\geq 1$ can also be rewritten in the following form:
$$q^b_n(t,x,y)=\int_0^t\int_{\R^d} p_0(t-r, x, z) (\S^b p_0)^{*, n}_z (r, z, y)\,dz\,dr.$$
Here $(\S^b p_0)^{*, n}_z (r, z, y)$ is defined in \eqref{e:s}.
Hence, for $T>t>s>t_0, x_1,x_2\in\R^d$ and $y\in\R^d,$ we have
$$\begin{aligned}
&|q^b_n(s,x_1,y)-q^b_n(t,x_2,y)|\\
\leq & \int_0^s\int_{\R^d}|p_0(s-r, x_1, z)-p_0(t-r, x_2, z)||(\S^b p_0)^{*, n}_z(r, z, y)|\,dz\,dr\\
&+\int_s^t\int_{\R^d} p_0(t-r, x_2, z)|(\S^b p_0)^{*, n}_z(r, z, y)|\,dz\,dr\\
=:& I+II.
\end{aligned}$$
It is known  that there are positive constants $c_1$ and $\theta$ so that
 for any $ t,s\in [t_0,T]$
and $x_i\in \R^d$ with $i=1,2,$
$$
|p_0(s,x_1,y)-p_0(t, x_2, y)|\leq
c_1\,  t_0^{-(d+\theta)/2}
\left(|t-s|^{1/2}+|x_1-x_2|\right)^\theta,
\quad y\in\R^d,
$$
we have by \eqref{e:3.9}, \eqref{e:f1} and Lemma \ref{n1}, for $\rho\in (0,s/2),$
\begin{equation}\label{e:h1}\begin{aligned}
I=&\int_0^{s-\rho}\int_{\R^d}|p_0(s-r, x_1, z)-p_0(t-r, x_2, z)||(\S^b p_0)^{*, n}_z(r, z, y)|\,dz\,dr\\
&+\int_{s-\rho}^s\int_{\R^d}|p_0(s-r, x_1, z)-p_0(t-r, x_2, z)||(\S^b p_0)^{*, n}_z(r, z, y)|\,dz\,dr\\
\leq &c_2 2^{-(n-1)}\rho^{-(d+\theta)/2}\left(|t-s|^{1/2}+|x_1-x_2|\right)^\theta\int_0^{s-\rho}\int_{\R^d} f_0(r,z,y)\,dz\,dr\\
&+c_2 2^{-(n-1)}(s-\rho)^{-(d+\beta)/2}\int_{s-\rho}^s\int_{\R^d}(p_0(s-r, x_1, z)+p_0(t-r, x_2, z))\,dz\,dr\\
\leq &c_3 2^{-(n-1)}\rho^{-(d+\theta)/2}\left(|t-s|^{1/2}+|x_1-x_2|\right)^\theta s^{1-\beta/2}+c_3 2^{-(n-1)}(s-\rho)^{-(d+\beta)/2}\rho.
\end{aligned}\end{equation}
Moreover, since $f_0(r,z,y)\leq s^{-(d+\beta)/2}$ for $r\in (s,t),$ we have
\begin{equation}\label{e:h2}
II \leq 2^{-(n-1)}\int_s^t\int_{\R^d} p_0(t-r, x_2, z)f_0(r,z,y)\,dz\,dr
 \leq 2^{-(n-1)}s^{-(d+\beta)/2}|t-s|.
 \end{equation}
Note that
$$
|q^b(s,x_1,y)-q^b(t,x_2,y)|
\leq |p_0(s,x_1,y)-p_0(t,x_2,y)|+\sum_{n=1}^\infty|q^b_n(s, x_1, y)-q^b_n(t, x_2, y)|.
$$
Then by taking $|t-s|$ and $|x_1-x_2|$ small, and then making $\rho$ small in \eqref{e:h1} and \eqref{e:h2}
 yields the conclusion of this Proposition.
\qed

\begin{prp}\label{H2}
For each  $0<t_0<T<\infty$ and $A>0$,
 the function $q^b(t,x,y)$ is uniform continuous in $y$ for
every $b$ with $\|b\|_\infty\leq A$ that satisfies \eqref{e:b}
 and for
all $(t,x)\in (t_0, T)\times \R^d.$
\end{prp}

\pf  By Theorem \ref{T:3.11} and the Chapman-Kolmogrov equation for $q^b(t,x,y)$,
it suffices to prove the theorem
for $A=A_0$ and $T=1$, where $A_0$ is the constant in Lemma \ref{u1}.

Define $P(s, x, y)=p_0(s,x)-p_0(s,y).$ For $s>0,$ we have
\begin{equation}\label{Lp}\begin{aligned}
&|\S^b p_0(s, y_1)-\S^b p_0(s, y_2)|\\
\leq &c_1\int_{\R^d}|P(s,y_1+h, y_2+h)-P(s, y_1, y_2)-\langle\nabla_{(y_1,y_2)} P(s, y_1, y_2), h
\1_{|h|\leq 1}\rangle|\frac{dh}{|h|^{d+\beta}}\\
\leq &c_1\int_{|h|\leq 1}|h|^2\sup_{\theta\in(0,1)}|\frac{\partial^2}{\partial y_1^2}p_0(s,y_1+\theta h)
-\frac{\partial^2}{\partial y_2^2}p_0(s,y_2+\theta h)|\frac{dh}{|h|^{d+\beta}}\\
&+c_1\int_{|h|>1} |p_0(s,y_1+h)-p_0(s, y_2+h)-p_0(s, y_1)+p_0(s,y_2)|\frac{dh}{|h|^{d+\beta}}\\
\leq &c_2 \sup_{y}|\frac{\partial^3}{\partial y^3}p_0(s,y)||y_1-y_2|\int_{|h|\leq 1}|h|^2\frac{dh}{|h|^{d+\beta}}
+c_2\sup_{y}|\frac{\partial}{\partial y}p_0(s,y)||y_1-y_2|\int_{|h|> 1}\frac{dh}{|h|^{d+\beta}}\\
\leq &c_3|y_1-y_2|[s^{-(d+3)/2}+s^{-(d+1)/2}],\\
\end{aligned}\end{equation}
where in the fourth inequality, we used $|\frac{\partial^3}{\partial y^3}p_0(s,y)|\leq c_3 s^{-(d+3)/2}$
which can be proved similarly by the argument in Lemma \ref{0}.

 Then for each $n\geq 1,$
 we have  by Lemma \ref{n1}, Lemma \ref{u1} and (\ref{Lp})
  that  for $(t, x, y)\in (t_0, 1)\times \R^d\times \R^d$ and  $\rho \in (0, t_0/2)$,
\begin{eqnarray*}
&&|q_n^b(t,x,y_1)-q_n^b(t,x,y_2)|\\
&\leq &\int_0^\rho\int_{\R^d}q_{n-1}^b(t-s,x,z)|\S_z^b p_0(s,z,y_1)-
\S_z^b p_0(s,z,y_2)|\,dz\,ds\\
&&+\int_\rho^t\int_{\R^d}q_{n-1}^b(t-s,x,z)|\S_z^b p_0(s,z,y_1)-\S_z^b p_0(s,z,y_2)|\,dz\,ds\\
&\leq &c_4 2^{-(n-1)}\int_0^\rho\int_{\R^d}  p_1
(t-s,x,z)|\S_z^b p_0(s,z,y_1)-\S_z^b p_0(s,z,y_2)|\,dz\,ds\\
&& +c_4 2^{-(n-1)} \int_\rho^t\int_{\R^d}
  p_1 (t-s,x,z)
\left| \S_z^b p_0(s,z-y_1)-\S_z^b p_0(s,z-y_2)\right|\,dz\,ds\\
&\leq &c_5 2^{-(n-1)}t_0^{-d/2}\int_0^\rho\int_{\R^d} \left(|\S_z^b p_0(s,z,y_1)|+|\S_z^b p_0(s,z,y_2)|\right)\,dz\,ds\\
&&+ c_5 2^{-(n-1)}\rho^{-(d+3)/2}|y_1-y_2|\int_\rho^t\int_{\R^d}
  p_1 (t-s,x,z)\,dz\,ds\\
&\leq & c_6 \, 2^{-(n-1)}\, t_0^{-d/2}\rho^{1-\beta/2}
+c_6 2^{-(n-1)}\rho^{-(d+3)/2}|y_1-y_2| .
\end{eqnarray*}
Therefore we have
$$\begin{aligned} &|q^b(t, x, y_1)-q^b(t, x, y_2)|\\
 \leq  &|p_0(t, x, y_1)-p_0(t, x, y_2)|+
 \sum_{n=1}^\infty c_62^{-(n-1)} \,   t_0^{-d/2}\rho^{1-\beta/2}
+\sum_{n=1}^\infty c_62^{-(n-1)}  \rho^{-(d+3)/2}|y_1-y_2|.
\end{aligned}$$
 By first taking $|y_1-y_2|$ small and then making $\rho$ small
 yields the desired uniform continuity of $q^b(t, x, y)$.
 \qed

\begin{thm}\label{T:m1}
Suppose $b$ is a bounded function on $\R^d\times\R^d$
satisfying \eqref{e:b} and \eqref{e:1.15}.
The kernel $q^b(t, x, y)$ uniquely determines a Feller process $X^b
=(X^b_t, t\geq 0, \P_x, x\in \R^d)$
on the canonical Skorokhod space
${\mathbb D}([0, \infty), \R^d)$
such that
$$ \E_x \left[ f(X^b_t)\right] =\int_{\R^d} q^b (t, x, y) f(y) dy
$$
for every bounded continuous function $f$ on $\R^d$.
The Feller process $X^b$ is conservative
 and has a L\'evy system $(J^b(x, y)dy, t)$,
where
$$ J^b(x, y)=\frac{  \, b(x, y-x)}{|x-y|^{d+\beta}}.
$$

\end{thm}
\pf When $b$ is a bounded function satisfying \eqref{e:b}, \eqref{e:5.1} and \eqref{e:1.15},
 the theorem has already been established via
Theorem \ref{T:1.2} and Propositions \ref{p2}-\ref{p4}.
We now remove the assumption \eqref{e:5.1}.
Suppose that $b(x, z)$ is a bounded
function that satisfies \eqref{e:b} and \eqref{e:1.15}.
Let  $\varphi$ be a non-negative smooth function with compact
support in $\R^d$ so that
$\int_{\R^d}\varphi (x) dx =1$. For each $n\geq 1$,
define $\varphi_n(x)=n^d  \varphi(nx)$ and
$$ k_n (x, z):= \int_{\R^d} \varphi_n (x-y) b(y, z) dy.
$$
Then $k_n$ is a function that satisfies \eqref{e:b}, \eqref{e:5.1} and \eqref{e:1.15}
with $\| k_n\|_\infty \leq \| b \|_\infty$.
By Theorems \ref{T:1.1}, \ref{T:1.2}, Propositions \ref{H1} and \ref{H2}, $q^{k_n}(t,x,y)$ is nonnegative, uniformly bounded
and equi-continuous on $[1/M, M]\times\R^d\times\R^d$ for each $M\geq 1$, then
there is a subsequence $\{n_j\}$ of $\{n\}$ so that $q^{k_{n_j}}(t,x,y)$ converges  boundedly and uniformly on compacts of $(0, \infty)\times \R^d\times \R^d$,
to some nonnegative continuous function $\overline q (t,x,y)$, which again satisfies \eqref{e:1.7}. Obviously, $\overline q(t,x,y)$ also satisfies the Chapman-Kolmogorov equation and $\int_{\R^d}\overline q(t,x,y)\,dy=1.$
By \eqref{e:3.21} and Theorem \ref{c1},
$$
q^{k_{n_j}}(t, x, y)
=p_0(t, x, y)+\int_0^t \int_{\R^d} q^{k_{n_j}} (t-s, x, z)
\S^{k_{n_j}}_z p_0(s, z, y)dz ds
$$
and
$$  q^{k_{n_j}}(t, x, y)\leq c_1\, p_1(t, c_2x, c_2y)
$$
for every $0<t\leq 1\wedge(A_0/\|b\|_\infty)^{2/(2 -\beta )}$
and $x, y\in \R^d$, where $c_k, k=1,2$ are positive constants that depend only
on $d, \beta$ and $\|b\|_\infty$. Letting $j\to \infty$, we have
by \eqref{L}, Lemma \ref{2'} and the dominated convergence theorem
that
$$
\overline q(t, x, y)=p_0 (t, x, y)+\int_0^t \int_{\R^d} \overline q (t-s, x, z)
\S^b_z p_0 (s, z, y)dy ds
$$
and $  \overline q(t, x, y)\leq c_1\, p_1(t, c_2x, c_2y)$
for every $0<t\leq 1\wedge(A_0/\|b\|_\infty)^{2/(2 -\beta )}$
and $x, y\in \R^d$. Hence we conclude from Theorem \ref{T:3.10}
that $\overline q(t, x, y)= q^b(t, x, y)$. This in particular implies
that $q^b(t, x, y)\geq 0$. So there is a Feller process $X^b$ having $q^b(t, x,y)$
as its transition density function.
The proof of Propositions \ref{p2}-\ref{p4}
only uses the condition (\ref{e:5.1}) through its implication
 that $q^b(t, x, y)\geq 0$. So in view of what we just established,
Propositions \ref{p2}-\ref{p4}  continue to
hold for $X^b$ under the current setting without the additional assumption \eqref{e:5.1}.
The proof of the theorem is now complete.
\qed

\medskip

For a Borel set $B\subset \R^d$, we define
 $\tau^b_B=\inf\{t>0: X^b_t\notin B\}$ and
  $\sigma^b_B:=\inf\{t\geq 0: X^b_t\in B\}.$

\begin{prp}\label{P:y1}
For each $A>0$ and $R_0>0$, there exists
a positive constant
$
\kappa =\kappa (d, \beta, A, R_0)< 32/9
$
so that
for every $b$ satisfying \eqref{e:b} and \eqref{e:1.15}
with $\|b\|_\infty \leq A,  x\in \R^d$ and $r\in (0,R_0),$
$$
\P_x \left(\tau^b_{B(x,r)}\leq \kappa r^{2} \right)
\leq \dfrac{1}{2} .
$$
\end{prp}

\pf  Let $f$ be a $C^2$ function taking values in $[0,
1]$ such that $f(0)=0$ and $f(u)=1$ if $|u|\geq 1.$ Set
$f_{x, r}(y)=f(\frac{y-x}{r})$.  Note that $f_{x, r}$ is a $C^2$
function taking values in $[0, 1]$ such that $f_{x,r}(x)=0$ and
$f_{x, r}(y)=1$ if $y\notin B(x, r)$. Moreover,
$$
\sup_{y\in\R^d}\left|\frac{\partial^{2}f_{x,r}(y)}{\partial y_i\partial y_j}\right| \leq r^{-2}\,
\sup_{y\in\R^d}\left|\frac{\partial^{2}f (y)}{\partial
y_i\partial y_j}\right|.
$$
Denote $\sum_{i, j=1}^d | \partial^2_{ij} f(x)|$ by
$|D^2 f(x)|$.
By Taylor's formula,
it follows that \begin{equation}\label{10}\begin{aligned} |\L^b
f_{x,r}(u)| &\leq |\Delta f_{x,r}(u)|+ c_1\int_{|h|\leq r} \left|f_{x,r}(u+h)-f_{x,r}(u)-\langle\nabla
f_{x,r}(u),h\rangle\right|\frac{dh}{|h|^{d+\beta}}\\
&\quad +\int_{|h|>r} |f_{x,r}(u+h)-f_{x,r}(u)|\frac{dh}{|h|^{d+\beta}}\\
&\leq r^{-2}|\Delta f(u)|+c_2\|D^2f\|_\infty r^{-2}\int_{|h|\leq 1} |h|^2\dfrac{dh}{|h|^{d+\beta}}
+c_2\|f\|_\infty \int_{|h|>r} \frac{dh}{|h|^{d+\beta}}\\
&\leq  c_3 (r^{-2}+r^{-\beta})\leq c_4r^{-2}, \quad r\in (0,R_0),
\end{aligned}\end{equation}
where $c_4=c_4(d,\beta, A, R_0)$ is a positive constant dependent on $R_0$.
Therefore, for each $t>0,$
$$\begin{aligned}
\P_x(\tau^b_{B(x,r)}\leq t)
&\leq \E_x \left[f_{x,r}(X^b_{\tau^b_{B(x,r)}\wedge t})\right]-f_{x,r}(x)\\
& =  \E_x \left[ \int_0^{\tau^b_{B(x,r)}\wedge t} \L^b f_{x,r}(X^b_s)\,ds \right]
\leq   c_4 \dfrac{t}{r^{2}}.
\end{aligned}$$
Set $\kappa=32/9\wedge (2c_4)^{-1},$ then
$$\P_x(\tau^b_{B(x,r)}\leq \kappa r^{2})\leq \frac{1}{2}.$$
\qed

Recall that $m_b=\inf_x {\rm essinf}_z b(x,z).$

\begin{prp}\label{x1}
For every $A>0$ and $R_0>0$,
there exists a constant $C_{18}=C_{18}(d,\beta, A, R_0)>0$
so that for every  $b$ satisfying \eqref{e:b} and \eqref{e:1.15}
with $\| b\|_\infty \leq A$,
$r\in (0, R_0]$ and $x, y\in \R^d$ with $|x-y|\geq 3r$,
$$
\P_x \big(\sigma^b_{B(y, r)}<\kappa r^{2}\big)
\geq C_{18}\, r^{d+2} \,
\frac{m_{b}}{  |x-y|^{d+\beta}}.
$$
\end{prp}

\pf  By Proposition \ref{P:y1} ,
$$
\E_x \left[ \kappa r^{2}\wedge \tau^b_{B(x,r)} \right]
\geq \kappa r^{2}\, \P_x \big (\tau^b_{B(x,r)}\geq
\kappa r^{2} \big)\geq \frac{1}{2}\kappa r^{2}.
$$
Thus by Proposition \ref{p4}, we have for $|x-y|\geq 3r,$
\begin{eqnarray*}
\P_x(\sigma^b_{B(y, r)}<\kappa r^{2})
  &\geq&
\P_x(X^b_{\kappa r^{2}\wedge \tau^b_{B(x,r)}}\in B(y,r))\\
&=& \E_x\int_0^{\kappa r^{2}\wedge\tau^b_{B(x,r)}}\int_{B(y,r)}
J^b(X^b_s, u)\,du\,ds\\
&\geq & c_1   \E_x \big[\kappa r^{2}\wedge
\tau^b_{B(x,r)}\big]\int_{B(y,r)} \frac{m_b}{  |x-y|^{d+\beta}}  du\\
&\geq&  c_2 \kappa r^{d+2}
\frac{m_{b}}{  |x-y|^{d+\beta}}.
\end{eqnarray*}
\qed

\begin{prp}\label{x1'}
For every $A>0, \lambda>0$ and $\ee>0$,
there exists a constant $C_{19}=C_{19}(d,\beta, A, \ee, \lambda)>0$
so that for every bounded $b$ that  satisfies \eqref{e:b}, \eqref{e:1.15} and \eqref{e:1.8}
with $\| b\|_\infty \leq A$,
 and  $3r\leq |x-y|\leq \lambda/3$,
$$
\P_x \big (\sigma^b_{B(y, r)}<\kappa r^{2}\big)
\geq C_{19}\, \frac{r^{d+2}}{|x-y|^{d+\beta}}.
$$
\end{prp}

\pf
By Propositions \ref{p4} and \ref{P:y1},
we have for $3r\leq |x-y|\leq \lambda/3$,
$$\begin{aligned}
\P_x \left(\sigma^b_{B(y, r)}<\kappa r^{2} \right)
&\geq \P_x \left( X^b_{\kappa r^{2}\wedge \tau^b_{B(x,r)}}
 \in B(y,r) \right)\\
&=\E_x\int_0^{\kappa r^{2}\wedge\tau^b_{B(x,r)}}\int_{B(y,r)}
J^b(X^b_s, u)\,du\,ds\\
&\geq c_1\E_x \left[ \kappa r^{2}\wedge
\tau^b_{B(x,r)} \right] \int_{B(y,r)}
\frac{\ee}{|x-y|^{d+\beta}} du\\
&\geq c_2\kappa r^{d+2} \, \frac{\ee}{|x-y|^{d+\beta}},
\end{aligned}$$
where the second inequality holds due to \eqref{e:1.8} and $|X^b_s-u|\leq 3 |x-y|\leq \lambda$ for $u\in B(y,r)$ and $X_s^b\in B(x,r)$.
\qed

\begin{thm}\label{t2'}
For every $A>0$ and
 any  $b$ satisfying  \eqref{e:b} and \eqref{e:1.15}  with $\| b\|_\infty\leq A$,
\eqref{e:1.14'} holds.
\end{thm}

\pf By \eqref{e:1.7}, we only need to prove the lower bound.
Let $\delta_0:=1\wedge (A_0/A)^{2/(2-\beta)}.$
\eqref{e:ql1} together with (\ref{e:1.5})   yields that for
any  $\|b\|_\infty\leq A$,
\begin{equation}\label{e:lq^b}
q^b(t,x,y)\geq c_0 t^{-d/2}  \quad \hbox{for }
t\in(0,\delta_0] \hbox{ and } \,|x-y|\leq 3t^{1/2},
\end{equation}
where $c_0=c_0(d, \beta)$ is a positive constant.
By \eqref{e:lq^b} and the usual chain argument, there are positive constants $c_1$ and $c_2$ so that
\begin{equation}\label{e:ql2}
q^b(t,x,y)\geq c_1p_0(t, c_2x, c_2y), \quad x,y\in\R^d, \,t\in (0, \delta_0].
\end{equation}

 For every $t\in (0, \delta_0]$,
by Proposition \ref{P:y1} and Proposition
\ref{x1} with  $R_0=1$,
$r=t^{1/2}/2$ and the strong Markov property of the process $X^b,$ we get for $|x-y|>3t^{1/2},$
\begin{eqnarray}\label{P}
&& \P_x(X^b_{2^{-2}\kappa t}\in B(y, t^{1/2})) \nonumber \\
&\geq& \P_x\left( X^b \hbox{ hits} B(y, t^{1/2}/2) \hbox{ before }
\frac{1}{4} \kappa t \hbox{ and stays there for at least }
 \frac{1}{4} \kappa t \hbox{ units of time} \right) \nonumber \\
&\geq & \P_x \left(\sigma^b_{B(y, t^{1/2}/2)}<\frac{1}{4}\kappa t \right) \inf_{z\in B(y, t^{1/2}/2)}
\P_z \left( \tau^b_{B(y, t^{1/2})}\geq \frac{1}{4}\kappa t\right) \nonumber \\
&\geq & \P_x \left( \sigma^b_{B(y, t^{1/2}/2)}<\frac{1}{4}\kappa t \right) \inf_{z\in B(y, t^{1/2}/2)}
\P_z \left( \tau^b_{B(z, t^{1/2}/2)}\geq \frac{1}{4}\kappa t \right) \nonumber \\
&\geq&  c_3\, t^{(d+2)/2}
\frac{m_b}{|x-y|^{d+\beta}}.
\end{eqnarray}
Here $c_3=c_3(d, \beta, A)$ is a positive constant.
 Hence,  by (\ref{e:lq^b}) and (\ref{P}), for $|x-y|>3t^{1/2}$ and $t\in (0,\delta_0],$
\begin{equation}\label{e:5.ql}\begin{aligned}
q^b(t,x,y)&\geq \int_{B(y,t^{1/2})} q^b(\frac{1}{4}\kappa t,x,z)q^b((1-\frac{1}{4}\kappa) t, z, y)\,dz\\
&\geq \inf_{z\in B(y,t^{1/2})}q^b((1-\frac{1}{4}\kappa)t,z,y)\P_x(X^b_{\frac{1}{4}\kappa t}\in B(y,t^{1/2}))\\
&\geq c_4t^{-d/2} \, t^{(d+2)/2}
 \frac{m_b}{|x-y|^{d+\beta}}\\
&\geq c_4
\frac{m_bt}{ |x-y|^{d+\beta}},
\end{aligned}\end{equation}
where $c_4=c_4(d,\beta, A)>0,$ the third inequality holds due to (\ref{e:lq^b}), (\ref{P}) and $|z-y|\leq t^{1/2}\leq 3((1-\frac{1}{4}\kappa)t)^{1/2}$
when $\kappa\leq 32/9$.
Finally, noting that $a\vee b\asymp a+b,$ \eqref{e:lq^b}, \eqref{e:ql2}, \eqref{e:5.ql}, \eqref{e:h} and the Chapman-Kolmogorov equation yields the desired lower
bound estimate. \qed

\begin{thm}\label{t2''}
For every $\lambda>0, \ee>0, A>0$ and  any  bounded $b$
satisfying \eqref{e:b}, \eqref{e:1.15} and \eqref{e:1.8} with $\| b\|_\infty\leq A$,
\eqref{e:1.14} holds.
\end{thm}

\pf Let $\delta_0:=1\wedge (A_0/A)^{2/(2-\beta)}\wedge (\lambda/9)^2.$
By  Theorem \ref{t2'} and Chapman-Kolmogorov equation, it suffices to prove there exist $c_k=c_k(d, \beta, A, \ee, \lambda)>0, k=1,2$
so that  $q^b(t,x,y)\geq c_1\overline p_\beta(t, c_2x, c_2y)$ for $t\in (0,\delta_0]$ and $|x-y|>3t^{1/2}$.

(i) First, we consider the case $\lambda/3 \geq |x-y|>3t^{1/2}.$
For every $t\in (0, \delta_0]$,
by Proposition
\ref{x1'} with
$r=t^{1/2}/2$ and the similar procedure in (\ref{P}),
\begin{equation}\label{P'}
\P_x \left( X^b_{\frac{1}{4}\kappa t}\in B(y, t^{1/2}) \right)
\geq c_3\, t^{(d+2)/2}\frac{1}{|x-y|^{d+\beta}}, \quad \lambda/3 \geq |x-y|>3t^{1/2} .
\end{equation}
Here $c_3=c_3(d,\beta, A, \ee, \lambda)$ is a positive constant.
 Hence,   by (\ref{e:lq^b}), (\ref{P'}) and the similar argument in \eqref{e:5.ql}, we have
\begin{equation}\label{l'b1}
q^b(t,x,y)
\geq c_4 \frac{t}{|x-y|^{d+\beta}}, \quad \lambda/3 \geq |x-y|>3t^{1/2}
\end{equation}
where $c_4=c_4(d,\beta, A, \ee, \lambda)>0.$

(ii) Next, we consider the case $ |x-y|>\lambda/3.$
 The following proof is similar to \cite[Theorem 3.6]{CKK}. For the reader's  convenience, we
spell out the details here.

Take $C_\ast=(\lambda/3) ^{-1}.$
Let $R:=|x-y|$ and $c_+=C_\ast\vee \delta_0^{-1}.$ Let $l\geq 2$ be a positive integer such that
$c_+R\leq l\leq c_{+}R+1$ and let $x=x_0, x_1, \cdots, x_l=y$ be such that $|x_i-x_{i-1}|\asymp R/l \asymp 1/c_+$
for $i=1,\cdots, l-1.$ Since $t/l\leq C_\ast R/l\leq C_\ast/c_+\leq \delta_0$ and $R/l\leq 1/c_+\leq \lambda/3 ,$
we have by  \eqref{e:lq^b}, \eqref{e:ql2} and \eqref{l'b1},
\begin{equation}\label{l'b}
q^b(t/l, x_i, x_{i+1})\geq c_5(t/l)^{-d/2}\wedge \left(p_0(t/l, c_6x_i, c_6x_{i+1})+\frac{t/l}{(R/l)^{d+\beta}}\right)
\geq c_7\left((t/l)^{-d/2}\wedge (t/l)\right)\geq c_7t/l.
\end{equation}
Let $B_i=B(x_i,\lambda/6 ),$ by (\ref{l'b}),
\begin{equation}\label{l'b2}\begin{aligned}
q^b(t,x,y)&\geq \int_{B_{l-1}}\cdots\int_{B_{1}} q^b(t/l,x,z_1)\cdots q^b(t/l,z_{l-1}, y)\,dz_1\cdots dz_{l-1}\\
&\geq (c_7t/l)^l\geq (c_8t/R)^{c_+R+1}\geq c_9(t/R)^{c_{10}R}\\
&\geq c_9\left(\frac{t}{|x-y|}\right)^{c_{10}|x-y|}.
\end{aligned}\end{equation}
By \eqref{l'b1}, \eqref{l'b2} and  the estimates of $\overline{p}_\beta$ in \eqref{e:trun1}-\eqref{e:trun2},
we get the desired conclusion.
\qed

 \bigskip
 \noindent{\bf Proof of Theorem \ref{T:1.3}.}
Theorem \ref{T:1.3} now follows from Theorems \ref{T:m1}, \ref{t2'} and \ref{t2''}.
\qed

\medskip

In the remainder of this section, we prove the result of Theorem \ref{T:1.4}, where $\L^b$ is a perturbation of $\Delta$
by finite range nonlocal operator $\S^b.$

\begin{prp}\label{y1}
For each $M>1,$ there exists
a positive constant $C_{20}=C_{20}(d, \beta, M, \lambda)$ so that
for every $b$ satisfying \eqref{e:b}, \eqref{e:1.15} and \eqref{e:1.10} and $r>0$,
$$\P_x(\tau^b_{B(x,r)}\leq t)\leq C_{20}\frac{t}{r^2}.$$
\end{prp}

\pf  Let $f$ be a $C^2$ function taking values in $[0,
1]$ such that $f(0)=0$ and $f(u)=1$ if $|u|\geq 1.$ Set
$f_{x, r}(y)=f(\frac{y-x}{r})$.  Note that $f_{x, r}$ is a $C^2$
function taking values in $[0, 1]$ such that $f_{x,r}(x)=0$ and
$f_{x, r}(y)=1$ if $y\notin B(x, r)$. Moreover,
$$
\sup_{y\in\R^d}\left|\frac{\partial^{2}f_{x,r}(y)}{\partial y_i\partial y_j}\right| \leq r^{-2}\,
\sup_{y\in\R^d}\left|\frac{\partial^{2}f (y)}{\partial
y_i\partial y_j}\right|.
$$
Denote $\sum_{i, j=1}^d | \partial^2_{ij} f(x)|$ by
$|D^2 f(x)|$.
By the conditions \eqref{e:1.15} and \eqref{e:1.10} $\sup_x b(x,z)\leq M1_{|z|\leq \lambda}(z)$,
it follows that \begin{equation}\label{10}\begin{aligned} |\L^b
f_{x,r}(u)| &\leq |\Delta f_{x,r}(u)|+ M\int_{|h|\leq \lambda} \left|f_{x,r}(u+h)-f_{x,r}(u)-\langle\nabla
f_{x,r}(u),h\rangle\right|\frac{dh}{|h|^{d+\beta}}\\
&\leq r^{-2}\|D^2f\|_\infty+M\|D^2f\|_\infty r^{-2}\int_{|h|\leq \lambda} |h|^2\dfrac{dh}{|h|^{d+\beta}}\\
&\leq  c r^{-2},
\end{aligned}\end{equation}
where $c=c(d,\beta, M, \lambda)$ is a positive constant independent of $r$.
Therefore, for each $t>0,$
$$\begin{aligned}
\P_x(\tau^b_{B(x,r)}\leq t)
&\leq \E_x \left[f_{x,r}(X^b_{\tau^b_{B(x,r)}\wedge t})\right]-f_{x,r}(x)\\
& =  \E_x \left[ \int_0^{\tau^b_{B(x,r)}\wedge t} \L^b f_{x,r}(X^b_s)\,ds \right]
\leq   c \dfrac{t}{r^{2}}, \quad r>0.
\end{aligned}$$
\qed

\begin{thm}\label{trun1}
For every $M>1$ and $\lambda>0$, there are positive constants
$C_k=C_k (d,\beta, M, \lambda), k=21, 22$
 such that
for any  $b$
with \eqref{e:b}, \eqref{e:1.15} and \eqref{e:1.10},
\begin{equation}\label{e:tr1}
q^b(t,x,y) \leq C_{21} \,t^{-d/2}\wedge [p_0(t, C_{22}x, C_{22}y)+ \overline p_{\beta} (t,C_{22}x,C_{22}y)], \quad t\in
(0,1],\, x,y\in\R^d.
\end{equation}
\end{thm}

\pf By  \eqref{e:3.26n}, \eqref{e:1.6t} and \eqref{e:trun1}, there are constants $c_k, k=1,\cdots, 4$ so that
\begin{equation}\label{e:trun4}
q^b(t,x,y)\leq c_1 p_1(t,c_2x,c_2y)\leq c_3 t^{-d/2}\wedge  [p_0(t,c_4x,c_4y)+\overline p_\beta (t, c_4x, c_4y)], \quad |x-y|\leq 1, \, t\in (0,1].
\end{equation}
 In the following, we will  prove that there exist $c_5$ and $c_6$ so that
\begin{equation}\label{e:trun3}
q^b(t,x,y)\leq c_5 \overline p_{\beta} (t, c_6x, c_6y), \quad |x-y|>1, \, t\in (0,1].
\end{equation}
In fact, by \eqref{e:trun2},
$\overline p_{\beta} (t, cx, cy)\leq c't^{-d/2}$ when $|x-y|>1$ and $t\in (0,1],$
so $\overline p_{\beta} (t, cx, cy) \asymp t^{-d.2}\wedge\overline p_{\beta} (t, cx, cy)$ in this case,
thus \eqref{e:trun3} implies that \eqref{e:tr1} holds for $|x-y|>1$ and $t\in (0,1].$

For each $\lambda>0,$ let $r>\lambda\vee 1$ be a constant to be chosen later.
First we will use induction method to prove that there is a constant $C_0$ so that
\begin{equation}\label{e:Hn}
\quad q^b(t,x,y)\leq C_0\left(\frac{t}{n}\right)^n, \quad \mbox{for} \, |x-y|\geq nr,\, t\in (0,1], \, n\geq 1.
\end{equation}
By  \eqref{e:3.26n} and \eqref{e:1.0}, we can find a constant $C_0$ so that
\begin{equation}\label{e:5.3}
q^b(t,x,y)\leq c_7p_1(t,c_8x,c_8y)\leq C_0\frac{t}{|x-y|^{d+\beta}} \leq C_0t, \quad |x-y|>1, \, t\in(0,1],
\end{equation}
where the second inequality holds since by \eqref{e:p00}, $p_0(t,x,y)\leq c\frac{t}{|x-y|^{d+2}}\leq c\frac{t}{|x-y|^{d+\beta}}$
for $|x-y|>1.$
Hence \eqref{e:Hn} naturally  holds for $n=1$.
Now fix $C_0$ and assume \eqref{e:Hn} holds for $n=m.$ Then $q^b(t,x,y)\leq C_0\left(t/m\right)^m$ for $|x-y|\geq mr$ and $t\in (0,1].$
Let $\tau$ be the first time that $X^b$ exits $B(x, r-\lambda)$ starting from $x$.
By the strong Markov property, for $|x-y|\geq (m+1)r,$
\begin{equation}\label{e:r0}
\begin{aligned}
q^b(t,x,y)&=\E_x[q^b(t-\tau, X^b_{\tau}, y); \tau<t]\\
&=\E_x\left[q^b(t-\tau, X^b_{\tau}, y); \tau<\frac{t}{m+1}\right]\\
&\quad +\E_x\left[q^b(t-\tau, X^b_{\tau}, y); \frac{t}{m+1}\leq \tau<t\right]\\
&=I+II.
\end{aligned}
\end{equation}

Noting that the jump of $X^b$ is not larger than $\lambda$ by L\'evy system formula Proposition \ref{p4} and the condition \eqref{e:1.10},
so starting from $x,$ $X^b_{\tau}$ will lie in $B(x,r),$ which implies the distance between $X^b_{\tau}$ and $y$ is bigger than $mr$
for $|x-y|>(m+1)r.$
Hence, by Proposition \ref{y1} and our assumption for $n=m$,
\begin{equation}\label{e:r1}
\begin{aligned}
I&\leq \P_x\left(\tau<\frac{t}{m+1}\right)\sup_{z\in B(x,r), s\leq t} q^b(s,z,y)\\
&\leq C_{20}\frac{t}{(m+1)(r-\lambda)^2}\cdot C_0(t/m)^m\\
&= \frac{C_0C_{20}}{(r-\lambda)^2}\left(\frac{m+1}{m}\right)^m \left(\frac{t}{m+1}\right)^{m+1}\\
&\leq \frac{C_0C_{20}K_1}{(r-\lambda)^2}\left(\frac{t}{m+1}\right)^{m+1}, m\geq 1,\, |x-y|>(m+1)r
\end{aligned}
\end{equation}
where $K_1>0$ is a positive constant independent of $m$ and $C_{20}$ is the constant in Proposition \ref{y1}.
On the other hand,  by Proposition \ref{y1} and the assumption for $n=m$,
\begin{equation}\label{e:r2}
\begin{aligned}
II&=\sum_{k=1}^m \E_x\left[q^b(t-\tau, X^b_{\tau}, y); \frac{kt}{m+1}\leq \tau<\frac{(k+1)t}{m+1}\right]\\
&\leq \sum_{k=1}^m \sup_{z\in B(x,r), s\leq (m+1-k)t/(m+1)} q^b(s,z,y)\cdot\P_x\left(\frac{kt}{m+1}\leq \tau<\frac{(k+1)t}{m+1}\right)\\
&\leq C_0\sum_{k=1}^m \left(\frac{1}{m}\frac{(m+1-k)t}{m+1}\right)^m\P_x\left(\tau<\frac{(k+1)t}{m+1}\right)\\
&\leq C_0C_{20}\sum_{k=1}^m \left(\frac{t}{m+1}\right)^m\left(\frac{m+1-k}{m}\right)^m\frac{(k+1)t}{(m+1)(r-\lambda)^2}\\
&=\frac{C_0C_{20}}{(r-\lambda)^2}\left(\frac{t}{m+1}\right)^{m+1}\sum_{k=1}^m (k+1)\left(\frac{m+1-k}{m}\right)^m\\
&=\frac{C_0C_{20}}{(r-\lambda)^2}\left(\frac{t}{m+1}\right)^{m+1}\sum_{k=1}^m (m+2-k)\left(\frac{k}{m}\right)^m, \quad |x-y|>(m+1)r.
\end{aligned}
\end{equation}
Define
$$S_m:=\sum_{k=1}^m (m+2-k)\left(\frac{k}{m}\right)^m, \quad m\geq 1.$$
For each $m\geq 1,$ define $f_m(u)=(m+2-u)(u/m)^m, u\in\R.$ Noting that
$f_m'(u)=\frac{u^{m-1}}{m^m}[m(m+2)-(m+1)u]>0$ for $u\in (0,m].$
Hence,
$$\begin{aligned}
S_m&=2+\sum_{k=1}^{m-1}(m+2-k)\left(\frac{k}{m}\right)^m\\
&\leq 2+\sum_{k=1}^{m-1}\int_k^{k+1}  (m+2-u)\left(\frac{u}{m}\right)^m\,du\\
&\leq 2+\frac{1}{m^m}\int_0^m (m+2-u)u^m\,du\\
&=2+\frac{1}{m^m}\left[\frac{(m+2)m^{m+1}}{m+1}-\frac{m^{m+2}}{m+2}\right]\\
&=2+\frac{m(3m+4)}{(m+1)(m+2)}\leq K_2, \quad m\geq 1,
\end{aligned}$$
where $K_2$ is a positive constant.
Then by \eqref{e:r2} and the above inequality, we have
$$II\leq \frac{C_0C_{20}}{(r-\lambda)^2}\left(\frac{t}{m+1}\right)^{m+1}S_m\leq \frac{C_0C_{20}K_2}{(r-\lambda)^2}\left(\frac{t}{m+1}\right)^{m+1}.$$
Combining this inequality with \eqref{e:r0} and \eqref{e:r1}, we have for $|x-y|\geq (m+1)r,$
$$q^b(t,x,y)=I+II\leq \frac{C_0C_{20}(K_1+K_2)}{(r-\lambda)^2}\left(\frac{t}{m+1}\right)^{m+1}.$$
Noting that  the constant $C_{20}$ in Proposition \ref{y1} is independent of $r,$ so we can take $r$ large enough
so that $\frac{C_{20}(K_1+K_2)}{(r-\lambda)^2}\leq 1.$ Thus, \eqref{e:Hn} holds for $n=m+1$ and thus for all $n\geq 1.$

When $|x-y|>2r,$ there exists $n$ so that $nr\leq|x-y|<(n+1)r,$
then by \eqref{e:Hn},
$$
q^b(t,x,y)\leq C_0\left(\frac{t}{n}\right)^n\\
\leq C_0\left(\frac{n+1}{n}\frac{rt}{|x-y|}\right)^{\frac{n}{n+1}\frac{|x-y|}{r}}
\leq C_0\left(\frac{2rt}{|x-y|}\right)^{\frac{|x-y|}{2r}}, \quad t\in (0,1].
$$
On the other hand, if $1<|x-y|\leq 2r,$ \eqref{e:5.3} shows that
$$q^b(t,x,y)\leq C_0 \frac{t}{|x-y|^{d+\beta}}\leq C_0 \frac{t}{|x-y|}\leq C_0 \left(\frac{2rt}{|x-y|}\right)^{\frac{|x-y|}{2r}}, t\in (0,1].$$
Hence,
\begin{equation}\label{e:100}
q^b(t,x,y)\leq C_0 \left(\frac{2rt}{|x-y|}\right)^{\frac{|x-y|}{2r}},\, |x-y|>1, \, t\in (0,1].
\end{equation}
Comparing \eqref{e:100} with \eqref{e:trun2}, we get \eqref{e:trun3} and the proof is complete.
\qed

 \bigskip
Theorem \ref{T:1.4} now follows from Theorem \ref{trun1}.

\bigskip

  {\bf Acknowledgement.} The author is grateful to Professor Z.-Q. Chen for his helpful comments.




\end{document}